\documentclass[a4paper,oneside,11pt]{article}

\usepackage{amsfonts,color}

\usepackage[final]{pdfpages}
\usepackage{graphicx, soul}
\usepackage{caption}
\usepackage{multirow}
\usepackage{hhline}
\usepackage{hyperref}

\newtheorem{dfn}{Definition}[section]
\newtheorem{tw}[dfn]{Theorem}
\newtheorem{prop}[dfn]{Proposition}
\newtheorem{rem}[dfn]{Remark}
\newtheorem{ex}[dfn]{Example}

\newtheorem{cor}[dfn]{Corollary}
\usepackage{amssymb} 
\usepackage{amsmath}
\numberwithin{equation}{section}

\renewcommand{\theequation}{\thesection.\arabic{equation}}

 \global\long\def\sbr#1{\left[ #1\right] }
 \global\long\def\cbr#1{\left\{  #1\right\}  }
 
 \global\long\def\rbr#1{\left(#1\right)}

 \global\long\def\E{\mathbb{E}}
 \global\long\def\P{\mathbb{P}}
 \global\long\def\R{\mathbb{R}}

 \global\long\def\dd#1{\textnormal{d}#1}

 \global\long\def\ra{\rightarrow}
 
 \global\long\def\ns{\infty}

\usepackage{anysize}
\usepackage{array}
\usepackage{enumerate}

\author{Micha\l \ Barski  \\ \small  Faculty of Mathematics, Warsaw University, Poland\\
 \small{\it m.barski@mimuw.edu.pl} \bigskip \\
\\
Rafa\l \ \L ochowski\footnote{The research of R{\L} was partially supported by the National Science Centre (Poland) under the grant agreements no.  2019/35/B/ST1/04292 and 2022/47/B/ST1/02114.}
\\ \small 
Department of Mathematics and Mathematical Economics,\\ \small Warsaw School of Economics, Poland\\ \small{\it rlocho@sgh.waw.pl}}

\title{\bf Affine term structure models driven by independent L\'evy processes}

\begin{document}

\maketitle

\begin{abstract}
\renewcommand{\theequation}{\arabic{equation}}

We characterize affine term structure models of non-negative short rate $R$ which may be obtained as solutions of autonomous SDEs driven by independent, one-dimensional L\'evy martingales, that is equations of the form
\begin{gather}\label{equation abstract}
\dd R(t)=F(R(t))\dd t+\sum_{i=1}^{d}G_i(R(t-))\dd Z_i(t), \quad R(0)=r_0\geq 0,\quad t>0,
\end{gather}
with deterministic real functions $F,G_1,...,G_d$ and independent one-dimensional L\'evy martingales $Z_1,...,Z_d$. 
Using a general result on the form of the generators of affine term structure models due to Filipovi\'c \cite{FilipovicATS}, it is shown, under the assumption that the Laplace transforms of the driving noises are regularly varying, that all possible solutions $R$ of \eqref{equation abstract} may be obtained also as solutions of autonomous SDEs driven by independent stable processes with stability indices in the range $(1,2]$. The obtained models include in particular the $\alpha$-CIR model, introduced by Jiao et al. \cite{JiaoMaScotti}, which proved to be still simple yet more reliable than the classical CIR model. Results on heavy tails of $R$ and its limit distribution in terms of the stability indices are proven. Finally, results of numerical calibration  of the obtained models to the market term structure of interest rates are presented and compared with the CIR and $\alpha$-CIR models. 
\end{abstract}


\section{Introduction}

\emph{Affine property} of a Markov process is (roughly saying) the property that the logarithm of the characteristic function of its transition kernel $p_t(x, \cdot)$ is given as an affine transformation of the initial state $x$.
This property was fundamental in the study of \emph{continuous state branching processes with immigration} (CBI) by Kawazu and Watanabe \cite{KawazuWatanabe}; this and other attractive analytical properties motivated  
Filipovi\'c  to bring in the pioneering paper \cite{FilipovicATS} \emph{affine processes}, which constitute the class of conservative CBI processes, in the field of finance.  Affine processes are widely used in various areas of mathematical finance, they appear in term structure models, credit risk modelling and are applied within the stochastic volatility framework. Solid fundamentals of affine processes in finance were laid down by Filipovi\'c \cite{FilipovicATS} and by Duffie, Filipovi\'c and Schachermeyer \cite{DuffieFilipovicSchachermeyer}. 
The results obtained in these papers settled a reference point for further research and proved the usefulness and strength  
of the Markovian approach. Missing questions on regularity and existence of c\`adl\`ag versions were answered by Cuchiero, Filipovi\'c and Teichmann \cite{CuchieroFilipovicTeichmann} and  Cuchiero and Teichmann \cite{CuchieroTeichmann}.
Dawson and Li \cite{DawsonLi06} gave a construction of CBI processes as strong solutions of systems of stochastic integral equations with random and non-Lipschitz coefficients, and jumps of Poisson type selected from some random sets. Such systems were further investigated by Fu and Li \cite{FuLi} and Dawson and Li \cite{DawsonLi12}.

The appearance of affine processes in finance has arguably started with the introduction of classical 
stochastic short rate models based on the Wiener process, like CIR (Cox, Ingersoll, Ross) \cite{CIR} and Vasi\v cek \cite{Vasicek}. Further research resulted in discovering new models, also with jumps; see, among others, Filipovi\'c \cite{FilipovicATS}, Dai and Singleton \cite{DaiSingleton}, Duffie and G\^arleanu \cite{DuffieGarleanu},  Barndorff-Nielsen and Shephard \cite{Bandorff-NielsenShepard}, Keller-Ressel and Steiner \cite{KellerRessel}, Jiao, Ma and Scotti \cite{JiaoMaScotti}.  A model framework based on stochastic dynamics is of  particular interest as it allows constructing discretization schemes enabling  e.g. Monte Carlo simulations which are essential for pricing exotic, i.e. path-dependent, derivatives.  A treatment of simulating schemes for affine processes and pricing methods can be found in \cite{Alfonsi1}.  Stochastic equations allow also to identify the number of random sources in the model which is of some use by calibration and hedging. Before we introduce the  stochastic integral equations of Dawson and Li \cite{DawsonLi06} let us state the form of the generator of a conservative CBI-process  (satisfying $p_t\rbr{x, [0, +\ns)} = 1$,  $t, x \ge 0$). 
A conservative CBI-process has (under the existence of the first moments assumption) the generator of the form
\begin{align}\label{generatorCBI}
\mathcal{A}f(x)=&c x f^{\prime\prime}(x)+(\beta x+b)f^\prime(x) +\int_{(0,+\infty)}\Big(f(x+y)-f(x)\Big)m(\dd y) \\[1ex] \nonumber
&+\int_{(0,+\infty)}\Big(f(x+y)-f(x)-f^\prime(x)y\Big)x\mu(\dd y), \quad x\geq 0,
\end{align}
where $c, b\geq 0$, $\beta\in\mathbb{R}$ and $m(\dd y)$, $\mu(\dd y)$ are nonnegative Borel measures on $(0,+\infty)$ satisfying
\begin{gather}\label{warunki na iary CBI}
\int_{(0,+\infty)}(1\wedge y)m(\dd y)+\int_{(0,+\infty)}(y\wedge y^2)\mu(\dd y)<+\infty.
\end{gather}
If $B$ is a standard Brownian motion, $N_m\rbr{\dd s, \dd y}$ and $N_{\mu} \rbr{\dd s, \dd y, \dd u}$ are Poisson random measures with intensities $\dd s \, m(\dd y)$ and $\dd s \, \mu(\dd y) \, \dd u$ respectively, $B$, $N_m\rbr{\dd s, \dd y}$ and $N_{\mu} \rbr{\dd s, \dd y, \dd u}$ are independent, and $\tilde{N}_{\mu} \rbr{\dd s, \dd y, \dd u}$ denotes the compensated $N_{\mu} \rbr{\dd s, \dd y, \dd u}$ measure then (under some technical assumption on  $m(\dd y)$) the stochastic equation 
\begin{align}
X(t) = & X_0 + \sqrt{2c} \int_0^t \sqrt{X(s)} \dd B_s +  \int_0^t \rbr{\beta X(s) + b} \dd s+  \int_0^t \int_0^{+\ns} y {N}_{m} \rbr{\dd s, \dd y} \nonumber \\
& + \int_0^t \int_0^{+\ns} \int_0^{X(s-)} y \tilde{N}_{\mu} \rbr{\dd s, \dd y, \dd u} ,\quad t\geq 0, \label{DawsonLi}
 \end{align}
has a unique non-negative strong solution, which is a CBI process with the generator given by \eqref{generatorCBI}, see \cite[Sect. 5]{DawsonLi06}, \cite{FuLi} and \cite[Sect. 3]{DawsonLi12}.

In this paper we focus on recovering from the form of their generator those affine processes, which are given as solutions of less general stochastic equations, namely the SDEs which are driven by a multidimensional L\'evy process with independent coordinates. Specifically, we focus on the equation 
\begin{gather}\label{rownanie 1}
\dd R(t)=F(R(t-))\dd t+\sum_{i=1}^{d}G_i(R(t-))\dd Z_i(t), \quad R(0)=r_0,\quad t>0,
\end{gather}
where $r_0$ is a nonnegative constant, $F$, $\{G_i\}_{i=1,2,...,d}$ are deterministic functions and $\{Z_i\}_{i=1,2,...,d}$ are independent L\'evy processes and martingales. A solution $R(t), t\geq 0$, if nonnegative, will be identified here with the short rate process which defines the bank account process by
$B(t):=\exp\rbr{\int_{0}^{t}R(s)ds},$ $t\geq 0$.
Related to the savings account are zero coupon bonds. Their prices form a family of stochastic processes
$P(t,T), t\in [0,T]$, parametrized by their maturity times $T\geq 0$.  The price of a bond with maturity $T$ at time $T$ is equal to its
nominal value, typically assumed, also here, to be $1$, that is $P(T,T)=1$. The family of bond prices is supposed to have the {\it affine structure}
\begin{gather}\label{affine model}
P(t,T)=e^{-A(T-t)-B(T-t) R(t)}, \quad 0\leq t\leq T,
\end{gather}
for some smooth deterministic functions $A$, $B:[0,+\infty) \rightarrow \mathbb{R}$. Hence, the only source of randomness in the affine model \eqref{affine model} is the short rate process $R$ given by \eqref{rownanie 1}.  As the resulting market constituted by $(B(t), \{P(t,T)\}_{T\geq 0})$
should exclude arbitrage, the discounted bond prices
$$
\hat{P}(t,T):=B^{-1}(t)P(t,T)=e^{-\int_{0}^{t}R(s)ds-A(T-t)-B(T-t)R(t)}, \quad 0\leq t\leq T,
$$
are supposed to be local martingales for each $T\geq 0$. This requirement affects in fact our starting equation \eqref{rownanie 1}. Thus the functions $F$, $\{G_i\}_{i=1,...,d}$ and the noise $Z=(Z_1,...,Z_d)$ should be chosen such that $\eqref{rownanie 1}$ has a nonnegative solution for any $x\geq 0$ and such that, for some functions $A$, $B:[0,+\infty) \rightarrow \mathbb{R}$ and each $T\geq 0$, \ $\hat{P}(t,T)$ is a local martingale on $[0,T]$. If this is the case, \eqref{rownanie 1} will be called to {\it generate an affine model} or to be a {\it generating equation}, for short.

The description of all generating equations with one-dimensional noise is well known, see Section \ref{Low-dimensional generating equations} for a brief summary. This paper deals with \eqref{rownanie 1} in the case $d>1$. The multidimensional setting makes the description of generating equations more involved due to the fact that two apparently different generating equations may have solutions which are Markov processes with identical generators. For brevity, we will call such solutions 'identical' or 'the same solutions'. The resulting bond markets are then the same, so such equations can be viewed as equivalent. The main results of the paper, i.e.
Theorem \ref{TwNiez}, Corollary \ref{cor o postaci mu} and Proposition \ref{prop canonical representation} imply under mild assumptions (regularly varying Laplace transforms of the driving noises) that  any generating equation \eqref{rownanie 1} has the same solution as that of the following equation 
\begin{gather}\label{canonical representation rownania}
\dd R(t)=(a R(t-)+b)\dd t+\sum_{k=1}^{g} d_k^{1/\alpha_k} R(t-)^{1/\alpha_k} \dd Z^{\alpha_k}_k(t), 
\end{gather}
with some $1\leq g\leq d$ and parameters  $a\in\mathbb{R}$, $b\geq 0$, $d_k>0$, $k=1,2,...,g$,  driven by independent
stable processes  $\{Z_k^{\alpha_k}\}$ with indices $\{\alpha_k\}$ such that $2\geq \alpha_1>\alpha_2>...>\alpha_g>1$.  All generating equations having the same solutions as 
\eqref{canonical representation rownania} form a class which we denote by
\begin{gather}\label{podklasy dla rownann}
\mathbb{A}_g(a,b;\alpha_1,\alpha_2,...,\alpha_g; \eta_1, \eta_2, ...,\eta_g),
\end{gather}
where $\eta_1 := d_1/2$ if $\alpha_1 = 2$, $\eta_i:=\frac{\Gamma(2-\alpha_i)}{\alpha_i(\alpha_i-1)} d_i$, if $\alpha_i \in (1,2)$, $i=1,...,g$; $\Gamma(\cdot)$ is the Gamma function. We call \eqref{canonical representation rownania} a {\it canonical representation} of \eqref{podklasy dla rownann}.
By changing values of the parameters in \eqref{podklasy dla rownann} one can thus
split all generating equations into disjoint subfamilies with a tractable canonical representation for each of them.  This classification is conceptually similar to that of Dai and Singleton \cite{DaiSingleton} obtained for the multivariate Wiener case. 

The number and structure of generating equations from the class \eqref{podklasy dla rownann} depend on the noise dimension in \eqref{rownanie 1}. As one may expect, this class becomes larger as $d$ increases.
In Section \ref{section Generalized CIR equations on a plane} we determine all generating equations on a plane by formulating
specific conditions for $F, G$ and $Z_1,Z_2$ in \eqref{rownanie 1}. For $d=2$ the class $\mathbb{A}_1(a,b;\alpha_1; \eta_1)$ 
consists of a wide variety of generating equations while $\mathbb{A}_2(a,b;\alpha_1,\alpha_2; \eta_1,\eta_2)$ turns out to be a singleton.
The passage to the case $d=3$ makes, however, $\mathbb{A}_2(a,b;\alpha_1,\alpha_2; \eta_1,\eta_2)$
a non-singleton. This phenomenon is discussed in Section \ref{section Example in higher dimensions}.

A tractable form of canonical representations is supposed to be an advantage for applications. 
One finds in \eqref{canonical representation rownania} with $g=1,\alpha_1=2$ the classical CIR model equation and  \eqref{canonical representation rownania} with $g=1,\alpha_1 \in (1,2)$ is the equation of  the \emph{stable CIR model}, considered e.g. in \cite{LiMa}, \cite{BarskiZabczykCIR}. One may expect that additional stable noise components improve the model of the bond market. 
For $g=2, \alpha_1=2$ and $\alpha_2 \in(1, 2)$, \eqref{canonical representation rownania} 
becomes the equation of the \emph{$\alpha$-CIR model} studied in \cite{JiaoMaScotti} (it is in place to mention that in  \cite{JiaoMaScotti} the authors introduce also much wider class of models, which they call \emph{$\alpha$-CIR integral type processes}, they contain all models from the classes $\mathbb{A}_g$, $g=1,2,\ldots$).  It was shown in \cite{JiaoMaScotti} that empirical behavior of the European sovereign bond market is closer to that implied by the $\alpha$-CIR model than by the CIR model
due to the permanent overestimation of the short rates by the latter one. The $\alpha$-CIR model allows also reconciling low interest rates with large fluctuations related to the presence of jump part whose tail fatness is controlled by the parameter $\alpha_2$. Exact asymptotics of tails of the short rate in the stable CIR model was given in \cite[Proposition 3.1]{LiMa}. In this paper we prove that the 
tail fatness of the short rate in the models from the class $\mathbb{A}_g(a,b;\alpha_1,...,\alpha_g; \eta_1,...,\eta_g)$ is controlled by the parameter $\alpha_g$. We also show estimations for the $p$-th moments of $R$ with $p<\alpha_g$ and characterize the limit distribution 
of $R(t)$ as $t\rightarrow +\infty$.

In the last part of the paper we focus on the calibration of canonical representations to market data. Into account are taken 
the spot rates of European Central Bank implied by the $AAA$ - ranked bonds. We compute numerically the fitting error for \eqref{rownanie 1}  in the Python programming language with $g$ in the range from $1$ up to $5$. 
This illustrates, in particular, the influence of $g$ on the reduction of fitting error which is always less than in the CIR model. The freedom of choice of stability indices makes the canonical model curves more flexible, hence with shapes better adjusted to the market curves. We observed that the $\alpha$-CIR model outperform the CIR model in this regard that it reduced the fitting error at least by $30\%$ in about $45\%$ considered cases, and in more than $65\%$ cases  considered the fitting error decreased by more than $10\%$. Unfortunately, addition of more noises (consideration of models from the classes $\mathbb{A}_g$, $g\ge3$) did not reduce the fitting error considerably; however, let us notice that addition of sufficiently fat tailed noise may be desirable from the risk management point of view since the noise with the fattest tail controls the tail fatness of the short rate.

The structure of the paper is as follows. Section \ref{section Preliminaries} contains a preliminary characterization of generating equations, i.e. Proposition \ref{prop wstepny}, which is a version of the result from  \cite{FilipovicATS} characterizing the generator of a Markovian short rate. This leads to a precise formulation of the problem studied in the paper. Further we 
describe one dimensional generating equations and discuss the non-uniqueness of generating equations in the multidimensional case.
Sect. \ref{section Classification of generating equations} is concerned with the classification of generating equations. Sect. \ref{section Noise with independent coordinates} contains the main results of the paper. In Sect. \ref{tail_fatness} we discuss the fatness of the tails of $R$ from the class $\mathbb{A}_g(a,b;\alpha_1,\alpha_2,...,\alpha_g; \eta_1,...,\eta_g)$ as well as the limit distribution of the short rate. Sections \ref{section Generalized CIR equations on a plane} and \ref{section Example in higher dimensions} are devoted to generating equations on a plane and an example in the three-dimensional case, respectively. In Sect. \ref{section Applications} we discuss the calibration of canonical representations. In the Appendix we prove Proposition \ref{prop wstepny} and Theorem \ref{tw d=2 independent coord.} .

\section{Preliminaries}\label{section Preliminaries}
In this section we present a version of the result on generators of
Markovian affine processes \cite{FilipovicATS}, see Proposition \ref{prop wstepny}, which is used for a precise formulation of the problem considered in the paper. We explain the meaning of the projections of the noise
and show in Example \ref{ex different eq ident sol} two different generating equations having the same projections, hence identical solutions.
For illustrative purposes we keep referring to the one-dimensional case where the forms of generating equations are well known, see Section
\ref{Low-dimensional generating equations} below. For the sake of notational convenience we often use a scalar product notation $\langle\cdot,\cdot\rangle$ in $\mathbb{R}^d$ and write \eqref{rownanie 1} in the form
 \begin{gather}\label{rownanie 2}
\dd R(t)=F(R(t-))\dd t+\langle G(R(t-)),\dd Z(t)\rangle, \quad R(0)=r_0\geq 0, \qquad t>0,
\end{gather}
where $G:=(G_1,G_2,...,G_d):[0,+\infty)\longrightarrow\mathbb{R}^d$ and 
 $Z:=(Z_1,Z_2,...,Z_d)$  is a L\'evy process in $\mathbb{R}^d$.

\subsection{Laplace exponents of L\'evy processes}
Let $Z$ be an $\mathbb{R}^d$-valued L\'evy process  with the characteristic triplet $(a,Q,\nu(\dd y))$ meaning that the characteristic function of $Z_t$, $t \ge 0$, reads 
$$  \mathbb{E}\sbr{e^{{\mathrm i} \langle \lambda,X(t)\rangle}} = \exp \rbr{ t\rbr{{\mathrm i}\langle \lambda, a\rangle - \frac{1}{2}\langle Q\lambda,\lambda\rangle+ \int_{\R^d \setminus \cbr{0}} e^{{\mathrm i} \langle \lambda, y\rangle} - 1 - {\mathrm i} \langle \lambda, y\rangle {\mathbf 1}_{\{\mid y\mid \le1\}}   \dd y }}, \quad \lambda \in \R^d.$$
We consider the case when $Z$ is a martingale. Consequently, $$\int_{\R^d} (|y|\wedge |y|^2)  \nu(\dd y) < +\ns,$$ the characteristic triplet of $Z$ is 
\begin{gather}\label{chrakterystyki Z}
\left(-\int_{\{\mid y\mid>1\}}y \ \nu(\dd y), \ Q, \ \nu(\dd y)\right)
\end{gather}  and we have the decomposition
$$
Z(t)=W(t)+X(t), \qquad X(t):=\int_{0}^{t}\int_{\mathbb{R}^d}y \ \tilde{\pi}(\dd s,\dd y), \quad t\geq 0,
$$
where $\tilde{\pi}(\dd s,\dd y) = {\pi}(\dd s,\dd y) - \dd s \nu(\dd y)$ is the compensated jump measure of $Z$ and $W$ is a $d$-dimensional Wiener process independent from $X$. The martingale $X$ will be called the jump part of $Z$. Its Laplace exponent $J_{X}$,  defined by 
$
\mathbb{E}\sbr{e^{-\langle \lambda,X(t)\rangle}}=e^{tJ_{X}(\lambda)}, 
$
has the  following representation
\begin{equation} \label{Jdef}
J_X(\lambda)=\int_{\mathbb{R}^d}(e^{-\langle\lambda,y\rangle}-1+\langle\lambda,y\rangle)\nu(\dd y),
\end{equation}
and is finite for $\lambda\in\mathbb{R}^d$ satisfying
$$
\int_{\mid y\mid>1}e^{-\langle \lambda,y\rangle}\nu(\dd y)<+\infty.
$$
By the independence of $X$ and $W$ the Laplace exponent $J_Z$ of $Z$ equals
\begin{equation} \label{LaplaceZ}
J_Z(\lambda)={{\frac{1}{2}\langle Q\lambda,\lambda\rangle+J_X(\lambda)}}.
\end{equation}

By a \emph{canonical stable martingale with index $\alpha=2$} (or a \emph{canonical $2$-stable martingale}) we will mean a one-dimensional standard Brownian motion $B$. By a \emph{canonical stable martingale with index $\alpha\in(1,2)$} (or a \emph{canonical $\alpha$-stable martingale}) we will mean a~real  L\'evy martingale $Z^\alpha(t), t\geq 0$,  with no Wiener part and the L\'evy measure of the form
$
\nu(\dd v):=\frac{1}{v^{\alpha+1}}\mathbf{1}_{\{v>0\}} \dd v
$. 
The Laplace exponent of a {canonical stable martingale with index $\alpha\in(1,2]$} reads
$
J_{B}(\lambda) = c_2 \lambda^2 \text{ if } \alpha = 2
$
and 
\begin{align}\label{LAplace exp alfa stabilny}
J_{Z^{\alpha}}(\lambda)&=\int_{0}^{+\infty}\left(e^{-\lambda v}-1+\lambda v\right)\frac{1}{v^{\alpha+1}} \dd v =
c_{\alpha} \lambda^\alpha, \quad \lambda \geq 0, \quad 
\end{align}
if $\alpha \in (1,2)$ with 
\begin{gather}\label{c alpha} 
c_2 = \frac{1}{2}, \quad c_\alpha=\frac{\Gamma(2-\alpha)}{\alpha(\alpha-1)}, \quad \alpha \in (1,2),
\end{gather}
where $\Gamma$ stands for the Gamma function.  For $\alpha \in (1,2)$ using  \cite[Property 1.2.15]{Taqqu} we 
obtain the following tail asymptotics
\begin{equation} \label{tails_Zzz}
\P\rbr{Z^{\alpha}(t) > z} \sim \frac{-t c_{\alpha}}{\Gamma\rbr{1-\alpha}} \frac{1}{z^{\alpha}} = \frac{t}{\alpha  z^{\alpha}} \text{ as } z \ra +\ns.
\end{equation}

\subsubsection{Projections of the noise}

For equation \eqref{rownanie 2} we consider the {\it projections} of $Z$ along $G$ given by
\begin{gather}\label{projection of Z}
Z^{G(x)}(t):=\langle G(x), Z(t)\rangle, \qquad x,t\geq 0.
\end{gather}
As linear transformations of $Z$, the projections form a family of real L\'evy processes parametrized by $x\geq 0$. 
If $Z$ is a martingale, then $Z^{G(x)}$ is a L\'evy martingale for any $x\geq 0$.
By the identity
$
\mathbb{E}\sbr{e^{- \gamma\cdot Z^{G(x)}(t)}}=\mathbb{E}\sbr{e^{-\langle \gamma G(x), Z(t)\rangle}}$, $\gamma\in\mathbb{R}, 
$
and \eqref{LaplaceZ} the Laplace exponent of $Z^{G(x)}$ equals 
\begin{gather}\label{Laplace ZG do uproszczenia}
J_{Z^{G(x)}}(\gamma)=J_Z(\gamma G(x))=\frac{1}{2}\gamma^2\langle Q G(x),G(x)\rangle+\int_{\mid y\mid>0}\left(e^{-\gamma \langle G(x),y\rangle}-1+\gamma\langle G(x),y\rangle\right)\nu(\dd y).
\end{gather}
Using the L\'evy measure $\nu_{G(x)}(\dd v)$ of $Z^{G(x)}$, which is the  image 
of the L\'evy measure $\nu(dy)$ under the linear transformation $y\mapsto \langle G(x), y\rangle$ given by 
\begin{gather}\label{nu_G}
\nu_{G(x)}(A):=\nu \{y \in \R^d: \langle G(x),y\rangle\in A \} , \quad A\in\mathcal{B}(\mathbb{R})
\end{gather}
we obtain that
\begin{gather}\label{Laplace ZG}
J_{Z^{G(x)}}(\gamma)=\frac{1}{2}\gamma^2\langle Q G(x),G(x)\rangle+\int_{\mid v\mid>0}\left(e^{-\gamma v}-1+\gamma v \right)\nu_{G(x)}(\dd v).
\end{gather}
Thus the characteristic triplet of the projection $Z^{G(x)}$ has the form
\begin{gather}\label{charakterystyki rzutu}
\left(-\int_{\mid v\mid>1}y \ \nu_{G(x)}(\dd v), \ \langle Q G(x),G(x)\rangle, \ \nu_{G(x)}(\dd v)\mid_{v\neq 0}\right).
\end{gather}
Above we used the restriction $\nu_{G(x)}(\dd v)\mid_{v\neq 0}$ by cutting off zero which may be an atom of $\nu_{G(x)}(\dd v)$.

\subsection{Preliminary characterization of generating equations}

In Proposition \ref{prop wstepny} below we provide a preliminary characterization for \eqref{rownanie 2} to be a generating equation. 
Note that the independence of coordinates of $Z$ is not assumed here. The central role play here the noise projections \eqref{projection of Z}. The result is deduced from Theorem 5.3 in \cite{FilipovicATS}, where the generator of a general non-negative Markovian  short rate process for affine models was characterized.

\begin{prop}\label{prop wstepny} 
Let $Z$ be a L\'evy martingale with characteristic triplet  \eqref{chrakterystyki Z} and $Z^{G(x)}$ be its projection \eqref{projection of Z} with the Le\'vy measure $\nu_{G(x)}(\dd v)$ given by \eqref{nu_G}.
\begin{enumerate}[(A)] 
\item Equation \eqref{rownanie 1} generates an affine model if and only if the following conditions are satisfied:
\begin{enumerate}[a)]
\item For each $x \ge 0$ the support of $\nu_{G(x)}$ is contained in $[0, +\ns)$ which means that $Z^{G(x)}$ has positive jumps only, i.e.  for each $t\geq 0$, with probability one,
\begin{gather}\label{Z^G positive jumps}
\triangle Z^{G(x)}(t):=Z^{G(x)}(t)-Z^{G(x)}(t-)=\langle G(x), \triangle Z(t)\rangle\geq 0.
\end{gather}
\item The jump part of $Z^{G(0)}$ has finite variation, i.e.
\begin{gather}\label{nu G0 finite variation}
\int_{(0,+\infty)}v \ \nu_{G(0)}(\dd v)<+\infty.
\end{gather}
\item The characteristic triplet \eqref{charakterystyki rzutu} of $Z^{G(x)}$ is linear in $x$, i.e.
\begin{align}\label{mult. CIR condition}
\frac{1}{2}\langle Q G(x), G(x)\rangle&=cx, \quad x\geq 0,\\[1ex]\label{rozklad nu G(x)}
\nu_{G(x)}(\dd v)\mid_{(0,+\infty)}&=\nu_{G(0)}(\dd v)\mid_{(0,+\infty)}+x\mu(\dd v), \quad x\geq 0,
\end{align}
for some $c\geq 0$ and a measure  $\mu(\dd v) \ \text{on} \ (0,+\infty) \ \text{satisfying}$
\begin{gather}\label{war calkowe na mu}
\int_{(0,+\infty)}(v \wedge v^2)\mu(\dd v)<+\infty.
\end{gather}
\item The function $F$ is affine, i.e.
\begin{gather}\label{linear drift}
F(x)=ax+b, \ \text{where} \ a\in\mathbb{R}, \ b\geq\int_{(1,+\infty)}(v-1)\nu_{G(0)}(\dd v) .
\end{gather}
\end{enumerate}
\item Equation \eqref{rownanie 1} generates an affine model if and only if the generator of $R$ is given by
\begin{align}\label{generator R w tw}\nonumber
\mathcal{A}f(x)=cx f^{\prime\prime}(x)&+\Big[ax +b+\int_{(1,+\infty)}(1 -v)\{\nu_{G(0)}(\dd v)+x\mu(\dd v)\}\Big]f^{\prime}(x)\\[1ex]
&+\int_{(0,+\infty)}[f(x+v)-f(x)-f^{\prime}(x)(1\wedge v)]\{\nu_{G(0)}(\dd v)+x\mu(\dd v)\}.
\end{align}
for $f\in\mathcal{L}(\Lambda)\cup C_c^2(\mathbb{R}_{+})$, where 
$\mathcal{L}(\Lambda)$ is the linear hull of $\Lambda:=\{f_\lambda:=e^{-\lambda x}, \lambda\in(0,+\infty)\}$
and $C_c^2(\mathbb{R}_{+})$ stands for the set of twice continuously differentiable functions with compact support in $[0,+\infty)$. The constants $a,b,c$ and the measures $\nu_{G(0)}(\dd v), \mu(\dd v)$ are those from part (A).
\end{enumerate}
\end{prop}

The poof of Proposition \ref{prop wstepny} is postponed to Appendix.

Note that conditions  \eqref{mult. CIR condition}-\eqref{rozklad nu G(x)} describe the distributions of the noise projections. In the sequel we use an equivalent formulation of  \eqref{mult. CIR condition}-\eqref{rozklad nu G(x)}  involving the Laplace exponents of \eqref{projection of Z}. Taking into account \eqref{Laplace ZG} we obtain the following.

\begin{rem}\label{rem warunki w eksp. Laplacea}
	The conditions \eqref{mult. CIR condition} and \eqref{rozklad nu G(x)}  are equivalent to the following decomposition of the Laplace exponent of $Z^G$:
	\begin{gather}\label{war na exp Laplaca}
	J_{Z^{G(x)}}(b)=cb^2x+J_{\nu_{G(0)}}(b)+x J_{\mu}(b), \quad b,x\geq 0,
	\end{gather}
	where
	\begin{gather}\label{def Jmu i Jnu0}
	J_{\mu}(b):=\int_{0}^{+\infty}(e^{-bv}-1+bv)\mu(\dd v), \quad J_{\nu_{G(0)}}(b):=\int_{0}^{+\infty}(e^{-bv}-1+bv)\nu_{G(0)}(\dd v).
	\end{gather}
\end{rem}

\subsubsection{Problem formulation}

In virtue of part $(A)$ of Proposition \ref{prop wstepny} we see that the drift $F$ of a generating equation is an affine function while 
the function $G$ and the noise $Z$ must provide projections $Z^{G(x)}, x\geq 0$ with particular distributions. 
Their characteristic triplets are characterized by a constant 
$c\geq 0$ carrying information on the variance of the Wiener part and two measures 
$\nu_{G(0)}(\dd v)$, $\mu(\dd v)$ describing jumps. 
A pair $(G,Z)$ for which the projections $Z^{G(x)}$ satisfy \eqref{Z^G positive jumps}-\eqref{war calkowe na mu}
will be called {\it a generating pair}. Note that the concrete forms of the measures $\nu_{G(0)}(\dd v)$, $\mu(\dd v)$
are, however, not specified. As for $Z$ with independent coordinates of infinite variation necessarily $G(0)=0$, see Proposition \ref{rem o G(0)=0}, and, consequently, $\nu_{G(0)}(\dd v)$ vanishes, our goal is to determine the measure $\mu(\dd v)$ in this case.

Having the required form of $\mu(\dd v)$ at hand one knows the distributions of the noise projections $Z^{G(x)}$ and, by  part $(B)$ of Proposition \ref{prop wstepny}, also the generator of the solution of \eqref{rownanie 2}. The generating pairs $(G,Z)$ can not be, however, 
uniquely determined, except the one-dimensional case.  This issue is discussed in 
Section \ref{Low-dimensional generating equations} and Section \ref{section Non-uniqueness in the multidimensional case} below.
For this reason we construct canonical representations - generating equations with noise projections corresponding to a given form of the measure $\mu(\dd v)$.

\subsubsection{One-dimensional generating equations}\label{Low-dimensional generating equations}

Let us summarize known facts on generating equations in the case $d=1$. If  $Z=W$ is a Wiener process, the only generating equation is the classical CIR equation
\begin{gather}\label{CIR equation}
\dd R(t)=(aR(t)+b)\dd t+C\sqrt{R(t)}\dd W(t), 
\end{gather}
with $a\in\mathbb{R}$, $b,C\geq 0$, see \cite{CIR}. The 
case with a general one-dimensional L\'evy process $Z$ was studied in \cite{BarskiZabczykCIR}, \cite{BarskiZabczyk} and
\cite{BarskiZabczykArxiv} with the following conclusion. If the variation of $Z$ is infinite and $G \not\equiv 0$, then $Z$ must be an $\alpha$-stable process with index $\alpha\in(1,2]$, with either positive or negative jumps only,  and \eqref{rownanie 1} has the form
\begin{gather}\label{CIR equation geenralized}
\dd R(t)=(aR(t-)+b)\dd t+C\cdot R(t-)^{{1}/{\alpha}}\dd Z^{\alpha}(t),
\end{gather}
with $a\in\mathbb{R}, b\geq 0$ and $C$ such that it has the same sign as the jumps of $Z^\alpha$. Clearly, for $\alpha=2$ equation \eqref{CIR equation geenralized} becomes \eqref{CIR equation}. If $Z$ is of finite variation then the noise enters \eqref{rownanie 1} in the additive way, that is 
\begin{gather}\label{Vasicek equation geenralized}
\dd R(t)=(aR(t-)+b)\dd t+C \ \dd Z(t).
\end{gather}
Here $Z$ can be chosen as an arbitrary process with positive jumps, $a\in\mathbb{R}, C\geq 0$ and 
$$
b\geq C \int_{0}^{+\infty}y \ \nu(\dd y),
$$
where $\nu(\dd y)$ stands for the L\'evy measure of $Z$.  The variation of $Z$ is finite, so is the right side above.
Recall, \eqref{Vasicek equation geenralized} with $Z$ being replaced by a Wiener process is the well known Vasi\v cek equation, see \cite{Vasicek}. Then the short rate is a Gaussian process, hence it takes negative values with positive probability. 
This drawback is eliminated by the jump version of  the Vasi\v cek equation \eqref{Vasicek equation geenralized}, where the solution never falls below zero.

It follows that the triplet $(c,\nu_{G(0)}(\dd v),\mu(\dd v))$ from  Proposition \ref{prop wstepny} takes for the equations above the following forms
\begin{enumerate} [a)]
\item $c \geq 0, \ \nu_{G(0)}(\dd v)\equiv 0, \ \mu(\dd v)\equiv 0$; \\[1ex]
This case corresponds to the classical CIR equation \eqref{CIR equation} where $c=\frac{1}{2}C^2$.
\item $c=0, \ \nu_{G(0)}(\dd v)\equiv0, \ \mu(\dd v)- \text{$\alpha$-stable}, \ \alpha\in(1,2)$;\\[1ex]
In this case \eqref{rownanie 2} becomes the stable CIR equation with $\alpha$-stable noise \eqref{CIR equation geenralized}.
\item $c=0, \ \nu_{G(0)}(\dd v)- \text{any measure on $(0,+\infty)$ of finite variation}, \ \mu(\dd v)\equiv 0$;\\[1ex]
Here \eqref{rownanie 2} becomes the generalized Vasi\v cek equation \eqref{Vasicek equation geenralized}.
\end{enumerate}
Note the one to one correspondence between the triplets $(c,\nu_{G(0)}(\dd v),\mu(\dd v))$  and generating pairs $(G,Z)$ which holds up to multiplicative constants.

\subsubsection{Non-uniqueness in the multidimensional case}\label{section Non-uniqueness in the multidimensional case}

In the case $d>1$ one should not expect a 1-1 correspondence between the triplets $(c,\nu_{G(0)}(\dd v),\mu(\dd v))$ and the generating equations \eqref{rownanie 2}. The reason is that the distribution of the noise projections $Z^{G(x)}$
does not determine the pair $(G,Z)$ in a unique way. Our illustrating example below shows two different equations 
driven by L\'evy processes with independent coordinates which provide the same short rate $R$. Note that the components of the process $\bar{Z}$ are \emph{not stable}.

\begin{ex}\label{ex different eq ident sol}

Let us consider the following two equations 
\begin{align}\label{ex pierwsze rownanie}
\dd R(t)&=\langle G(R(t-)), \dd Z(t) \rangle, \quad R(0)=R_0,\quad t\geq 0,\\[1ex]\label{ex drugie rownanie}
d\bar{R}(t)&=\langle \bar{G}(\bar{R}(t-), \dd \bar{Z}(t))\rangle, \quad \bar{R}(0)=R_0,\quad t\geq 0,
\end{align}
where 
$$
G(x):=2^{-1/\alpha}\cdot (x^{1/ \alpha}, x^{1/ \alpha}), \quad Z:=(Z^\alpha_1, Z^\alpha_2), 
$$
and 
$$
\bar{G}(x):=(x^{1/ \alpha}, x^{1/ \alpha}), \quad \bar{Z}:=(\bar{Z}_1, \bar{Z}_2),
$$
with a fixed index $\alpha\in(1,2)$. We assume that $Z^{\alpha}_1, Z^{\alpha}_2$ are independent canonical stable martingales with index $\alpha$ while $\bar{Z}_1, \bar{Z}_2$ are independent martingales with L\'evy measures
\[
\nu_1(\dd v) =  \frac{\dd v}{v^{\alpha +1}}\mathbf{1}_E(v), \quad \nu_2(\dd v) = \frac{\dd v}{v^{\alpha +1}} \mathbf{1}_{[0, +\ns) \setminus E}(v),
\]
respectively, where $E$ is a Borel subset of $[0, +\ns)$ such that 
$$
|E| = \int_{E} \dd v >0, \quad  \text{and} \quad |[0, +\ns) \setminus E| = \int_{[0, +\ns) \setminus E} \dd v  >0.
$$ 
The projections related to \eqref{ex pierwsze rownanie} and \eqref{ex drugie rownanie} take the forms
\begin{align*}
Z^{G(x)}(t)&=\langle G(x), Z(t)\rangle=x^{1/\alpha} 2^{-1/\alpha}(Z^\alpha_1(t)+Z^{\alpha}_2(t)), \quad x,t\geq 0,\\
\bar{Z}^{\bar{G}(x)}(t)&=\langle \bar{G}(x), \bar{Z}(t)\rangle=x^{1/\alpha} (\bar{Z}_1(t)+\bar{Z}_2(t)),\quad x,t\geq 0.
\end{align*}
Since both processes $2^{-1/\alpha} (Z^\alpha_1+Z^{\alpha}_2)$ and $\bar{Z}_1+\bar{Z}_2$ are canonical stable martingales with index
$\alpha$ we obtain that $(G, Z)$ and $(\bar{G}, \bar{Z})$ are generating pairs with the same solutions. 

It follows, in particular, that the noise coordinates of a generating equation do not need to be stable processes.
\end{ex}

\section{Classification of generating equations}\label{section Classification of generating equations}
\subsection{Main results}\label{section Noise with independent coordinates}

This section deals with equation \eqref{rownanie 2} in the case when the coordinates of the martingale $Z$ are independent. In view of Proposition \ref{prop wstepny}  we are interested in characterizing possible distributions of projections $Z^G$ over all generating pairs $(G,Z)$.  By \eqref{Z^G positive jumps} the jumps of the projections are necessarily positive. As the coordinates of $Z$ are independent, they do not jump together. Consequently, we see that, for each $x\geq 0$ and $t \ge 0$
$$
\triangle Z^{G(x)}(t)=\langle G(x),\triangle Z(t)\rangle >0
$$
holds if and only if, for some $i=1,2,...,d$,
\begin{gather}\label{mucha}
G_i(x)\triangle Z_i(t)>0, \quad \triangle Z_j(t)=0, j\neq i.
\end{gather}
Condition \eqref{mucha} means that $G_i(x)$ and $\triangle Z_i(t)$ are of the same sign. We can consider only the case when both are positive, i.e.
$$
G_i(x)\geq 0, \quad i=1,2,...,d, \ x\geq 0, \qquad \triangle Z_i(t)\geq 0, \quad t> 0,
$$
because the opposite case can be turned into this one by replacing $(G_i,Z_i)$ with $(-G_i,-Z_i)$, $i=1,...,d$. The L\'evy measure $\nu_i(\dd y)$ of $Z_i$ is thus concentrated on $(0,+\infty)$ and, in view of \eqref{LaplaceZ}, the Laplace exponent of $Z_i$ takes the form
\begin{gather}\label{Laplace Zi}
J_i(b):=\frac{1}{2}q_{ii} b^2+\int_{0}^{+\infty}(e^{-b v}-1+b v)\nu_i(\dd v),\quad b \geq 0, \ i=1,2,...,d,
\end{gather}
with $q_{ii}\geq 0$. Recall, $q_{ii}$ stands on the diagonal of $Q$ - the covariance matrix of the Wiener part of $Z$.  
We will assume that $J_i, i=1,2,...,d$ are  {\it regularly varying at zero}. Recall, this means that
$$
\lim_{x\rightarrow 0^+}\frac{J_i(bx)}{J_i(x)}=\psi_i(b), \quad b> 0,\qquad i=1,2,...,d,
$$
for some function $\psi_i$. In fact $\psi_i$ needs to be a power function, i.e.
$$
\psi_i(b)=b^{\alpha_i}, \quad b>0,
$$
with some $-\infty< \alpha_i<+\infty$ and $J_i$ is called to vary regularly with index $\alpha_i$, see \cite{Bingham}.

The distribution of noise projections are described by the following result.

\begin{tw} \label{TwNiez} Let $Z_1,...,Z_d$ be independent coordinates of the L\'evy martingale $Z$ in $\R^d$. Assume that $Z_1,...,Z_{d}$  satisfy
\begin{equation} \label{ass1}
\triangle Z_i(t)\geq 0 \text{ a.s.  for } t>0  \text{ and } Z_i \ \text{is of infinite variation}
\end{equation}
or 
\begin{equation} \label{ass2}
\triangle Z_i(t)\geq 0\text{ a.s.  for } t>0 \text{ and } G(0)=0.
\end{equation}
Further, let us assume that for all $i=1,\ldots, d$ the Laplace exponent  \eqref{Laplace Zi} of $Z_i$ varies regularly at zero and the components of the function  $G$ satisfiy
$$
G_i(x)\geq 0, \ x\in[0,+\infty), \quad G_i \ \text{is continuous on } [0,+\infty).
$$
Then \eqref{rownanie 2} generates an affine model if and only 
 if $F(x)=ax+b$, $a\in\mathbb{R}, b\geq 0$, and the Laplace exponent $J_{Z^{G(x)}}$ of $Z^{G(x)}=\langle G(x), Z\rangle$ is of the form 
\begin{gather}\label{postac J_ZG przy niezaleznych}
J_{Z^{G(x)}}(b) = x \sum_{k=1}^g\eta_{k}b^{\alpha_{{k}}}, \quad  \eta_{k}> 0, \quad \alpha_k\in(1,2],  \quad k=1,2,\ldots,g,
\end{gather}
with some $1 \le g \le d$ and $\alpha_k\neq\alpha_j$ for $k\neq j$.
\end{tw}

Theorem \ref{TwNiez} allows determining the form of the measure $\mu(\dd v)$ in Proposition \ref{prop wstepny}.

\begin{cor}\label{cor o postaci mu} 
Let the assumptions of Theorem \ref{TwNiez} be satisfied. If equation \eqref{rownanie 2} generates an affine model 
then the function $J_\mu$ defined in \eqref{def Jmu i Jnu0} takes the form
\begin{gather}\label{J mu postaaaccccc}
J_{\mu}(b)  = \sum_{k=l}^{g}\eta_{k}b^{\alpha_{{k}}}, \quad l\in \{1,2\}, \quad  \eta_{k}> 0, \quad \alpha_k\in(1,2),  \quad k=l,l+1,\ldots,g,
\end{gather}
with $1 \le g \le d$, $2>\alpha_l>...>\alpha_g>1$ (for the case $l=2, g=1$ we set $J_{\mu}\equiv 0$, which means that $\mu(\dd v)$ disappears). Above $l=2$ if $\alpha_1=2$ and $l=1$ otherwise. This means that $\mu(\dd v)$ is a weighted sum of $g+1-l$ stable measures with indices $\alpha_l,...,\alpha_g\in(1,2)$, i.e.
\begin{gather}\label{miary w klasach}
\mu(\dd v)=\tilde{\mu}(\dd v):= \frac{d_l}{v^{1+\alpha_l}}\mathbf{1}_{\{v>0\}}\dd v+...+\frac{d_g}{v^{1+\alpha_g}}\mathbf{1}_{\{v>0\}} \dd v,
 \end{gather}
with $d_i=\eta_i/ c_{\alpha_i},i=l,...,g$, where $c_{\alpha_i}$ is given by \eqref{c alpha} .
\end{cor}

Note that each generating equation can be identified by the numbers $a,b$ appearing in the formula for the function $F$ and $\alpha_1,...,\alpha_g; \eta_1,...,\eta_g$ from \eqref{postac J_ZG przy niezaleznych}. Since $\nu_{G(0)}(\dd v)=0$, see  Proposition \ref{rem o G(0)=0} in the sequel, the related generator of $R$ takes, by \eqref{generator R w tw}, the form
\begin{align}\label{generator niezalezne}\nonumber
 \mathcal{A}f(x)=cx f^{\prime\prime}(x)&+\Big[x\Big(a+\int_{(1,+\infty)}(1 -v)x\tilde{\mu}(\dd v)\Big)+b\Big]f^{\prime}(x)\\[1ex]
 &+\int_{(0,+\infty)}[f(x+v)-f(x)-f^{\prime}(x)(1\wedge v)]x\tilde{\mu}(\dd v),
\end{align}
with $\tilde{\mu}$ in \eqref{miary w klasach}. Recall, the constant $c$ above comes from the condition 
\begin{gather}\label{FG w klasach}
\frac{1}{2}\langle Q G(x),G(x)\rangle=cx, \quad  \quad x\geq 0,
\end{gather}
and, in view of Remark \ref{rem warunki w eksp. Laplacea}, $c=\eta_1$ if $\alpha_1=2$ and $c=0$ otherwise.
The class of processes with generator of the form \eqref{generator niezalezne} is denoted as in Eq. \eqref{podklasy dla rownann} by $\mathbb{A}_g(a,b;\alpha_1,...,\alpha_g; \eta_1,...,\eta_g)$.

 Note that the existence of the process being the strong, unique solution of \eqref{DawsonLi} with $\mu$ given by \eqref{miary w klasach}, $m\equiv 0$ and the generator given by \eqref{generator niezalezne} is guaranteed by \cite[Theorem 3.1]{DawsonLi12}.

\begin{prop}[Canonical representation of $\mathbb{A}_g(a,b;\alpha_1,...,\alpha_g; \eta_1,...,\eta_g)$]\label{prop canonical representation}
Let $R$ be the solution of \eqref{rownanie 2} with $F$, $G$ and $Z$ satisfying the assumptions of Theorem \ref{TwNiez}. 
 Let $\tilde{Z}=(\tilde{Z}^{\alpha_1},\tilde{Z}^{\alpha_2},...,\tilde{Z}^{\alpha_g})$ be a L\'evy martingale with independent coordinates which are canonical stable martingales with indices $\alpha_k, k=1,2,...,g$, respectively, and $\tilde{G}(x)=(d^{1/\alpha_1}_1 x^{1/\alpha_1},...,d^{1/\alpha_g}_g x^{1/\alpha_g})$, $x \ge 0$, where $d_k:=\eta_k/c_{\alpha_k}$ and $c_{\alpha_k}$  are given by \eqref{c alpha}, $k=1,2,...,g$.
Then 
$$
J_{Z^{G(x)}}(b)=J_{\tilde{Z}^{\tilde{G}(x)}}(b), \quad b,x\geq 0.
$$
Consequently, if $\tilde{R}$ is the solution of the equation 
\begin{gather}\label{rownanie sklajane}
\dd \tilde{R}(t)=(a\tilde{R}(t-)+b) \dd t+\sum_{k=1}^{g}d_k^{1/{\alpha_k}} \tilde{R}(t-)^{1/{\alpha_k}}\dd \tilde{Z}^{\alpha_k}(t),
\end{gather}
then the generators of $R$ and $\tilde{R}$ are equal.
\end{prop}

Equation \eqref{rownanie sklajane} will be called the {\it canonical representation} of the class $\mathbb{A}_g(a,b;\alpha_1,...,\alpha_g; \eta_1,...,\eta_g)$. The existence and uniqueness of the strong solution of \eqref{rownanie sklajane} follows for example from \cite[Theorem 5.3]{FuLi}.

\vskip2ex
\noindent
{\bf Proof:} By \eqref{postac J_ZG przy niezaleznych} we need to show that 
\begin{gather*}\label{rozklad rzutu konstrukcja}
J_{\tilde{Z}^{\tilde{G}(x)}}(b)=x\sum_{k=1}^{g}\eta_k b^{\alpha_k}, \quad b,x\geq 0.
\end{gather*}
Recall, the Laplace exponent of $\tilde{Z}^{\alpha_k}_k$ equals
$J_k(b)=c_{\alpha_k}b^{\alpha_k}, k=1,2,...,g$. By independence and the form of $\tilde{G}$ we have
\begin{align*}
J_{\tilde{Z}^{\tilde{G}(x)}}(b)&=\sum_{k=1}^{g}J_k(b\tilde{G}_k(x))=\sum_{k=1}^{g}c_{\alpha_k} b^{\alpha_k}d_kx=x\sum_{k=1}^{g}\eta_k b^{\alpha_k}, \quad b,x\geq 0,
\end{align*}
as required. The second part of the thesis follows from Proposition \ref{prop wstepny}(B).\hfill$\square$

\vskip1ex
Clearly, in the case $d=1$ the noise dimension can not be reduced, so $g=d=1$ and $\mathbb{A}_1(a,b; 2;\eta_1)$ corresponds to the classical CIR equation \eqref{CIR equation} while $\mathbb{A}_1(a,b; \alpha;\eta_1), \alpha\in(1,2)$ to its generalized version 
\eqref{CIR equation geenralized}. Both classes are singletons and \eqref{CIR equation}, \eqref{CIR equation geenralized} are their canonical representations. The $\alpha$-CIR equation from \cite{JiaoMaScotti} is a canonical representation of the class $\mathbb{A}_2(a,b; 2,\alpha; \eta_1,\eta_2)$ with $\alpha\in(1,2)$.

\subsubsection{Proofs}
The proofs of Theorem \ref{TwNiez} and Corollary \ref{cor o postaci mu} are preceded by two auxiliary results, i.e. Proposition \ref{bounds_alpha} and
Proposition \ref{rem o G(0)=0}. The first one provides some useful estimation for the function 
\begin{gather}\label{J}
J_{\rho}(b):=\int_{0}^{+\infty}(e^{-bv}-1+bv)\rho(\dd v), \quad b\geq 0,
\end{gather}
where the measure $\rho(\dd v)$ on $(0,+\ns)$ satisfies
\begin{gather}\label{nuJ}
0 < \int_0^{+\ns} \rbr{v^2\wedge v} \rho\rbr{\dd v} < +\ns.
\end{gather}
The second result shows that if all components of $Z$ are of infinite variation then $G(0)=0$.

\begin{prop} \label{bounds_alpha}
Let $J_{\rho}$ be a function given by \eqref{J} where the measure  $\rho$ satisfies \eqref{nuJ}. Then the function
$
(0,+\ns) \ni b \mapsto {J_{\rho}(b)}/{b}$ is strictly increasing and $\lim_{b \ra 0+}J_{\rho}(b)/b = 0$, while the function $(0,+\ns) \ni b \mapsto {J_{\rho}(b)}/{b^2}$
is strictly decreasing and $\lim_{b \ra +\ns}J_{\rho}(b)/b^2 = 0$. This yields, in particular, that, for any $b_0 >0$, 
\begin{gather}\label{oszacowania dwustronne J}
\frac{J_{\rho}\rbr{b_0}}{b_0^2}b^2 < J_{\rho}(b) < \frac{J_{\rho}\rbr{b_0}}{b_0}b, \quad b\in \rbr{0, b_0}.
\end{gather}
\end{prop} 
\noindent
{\bf Proof:} Let us start from the observation that the function 
$$
t \mapsto \frac{(1-e^{-t})t}{e^{-t}-1+t}, \quad t\geq 0,
$$
is strictly decreasing, with limit $2$ at zero and $1$ at infinity. This implies
 \begin{equation} \label{oszH}
 (e^{-t}-1+t) < (1-e^{-t})t < 2 (e^{-t}-1+t), \quad t \in (0, +\ns),
 \end{equation}
and, consequently,
$$
\int_{0}^{+\infty}(e^{-bv}-1+bv)\rho(\dd v) < \int_{0}^{+\infty}(1-e^{-bv})bv\ \rho(\dd v) < 2\int_{0}^{+\infty}(e^{-bv}-1+bv)\rho(\dd v), \quad b >0.
$$
This means, however, that
$$
J_{\rho}(b) < bJ_{\rho}^\prime(b) < 2J_{\rho}(b), \quad b > 0.
$$
So, we have
$$
\frac{1}{b} < \frac{J_{\rho}^\prime(b)}{J_{\rho}(b)}=\frac{d}{db}\ln J_{\rho}(b) < \frac{2}{b}, \quad b>0,
$$
and integration over some interval $[b_1,b_2]$, where $b_2 > b_1>0$, yields
$$
\ln b_2 - \ln b_1  < \ln J_{\rho}\rbr{b_2}-\ln J_{\rho}\rbr{b_1} < 2 \ln b_2 - 2 \ln b_1
$$
which gives that 
$$
\frac{J_{\rho}\rbr{b_2}}{b_2} > \frac{J_{\rho}\rbr{b_1}}{b_1}, \quad \frac{J_{\rho}\rbr{b_2}}{b_2^2} < \frac{J_{\rho}\rbr{b_1}}{b_1^2}.
$$

To see that $\lim_{b \ra 0+} {J_{\rho}\rbr{b}}/{b} = 0$ it is sufficient to use de l'H\^opital's rule,  \eqref{nuJ} and dominated convergence
$$
\lim_{b \ra 0+} \frac{J_{\rho}\rbr{b}}{b} = \lim_{b \ra 0+} {J'_{\rho}\rbr{b}} = \lim_{b \ra 0+}  \int_{0}^{+\infty}(1-e^{-bv}) v\ \rho(\dd v) = 0.
$$

To see that $\lim_{b \ra +\ns} {J_{\rho}\rbr{b}}/{b^2} = 0$ we also use de l'H\^opital's rule,  \eqref{nuJ} and dominated convergence.
If $\int_0^{+\ns} v\ \rho\rbr{\dd v} < +\ns$, then we have
$$
\lim_{b \ra +\ns} \frac{J_{\rho}\rbr{b}}{b^2} = \lim_{b \ra +\ns} \frac{J_{\rho}'\rbr{b}}{2b} = \frac{\int_0^{+\ns} v \rho\rbr{\dd v}}{+\ns}= 0.
$$
If $\int_0^{+\ns} v\ \rho\rbr{\dd v} = +\ns$ then we apply de l'H\^opital's rule twice and obtain
$$
\lim_{b \ra +\ns} \frac{J_{\rho}\rbr{b}}{b^2} = \lim_{b \ra +\ns} \frac{J_{\rho}'\rbr{b}}{2b} = \lim_{b \ra +\ns} \frac{J_{\rho}''\rbr{b}}{2} = \frac{1}{2} \lim_{b \ra +\ns} \int_{0}^{+\infty}e^{-bv} v^2\ \rho(\dd v) = 0.
$$
\hfill $\square$

\begin{prop}\label{rem o G(0)=0}
If $(G,Z)$ is a generating pair and all components of $Z$ are of infinite variation then $G(0)=0$. 
\end{prop}
\noindent 
{\bf Proof:}  Let $(G,Z)$ be a generating pair. Since the components of $Z$ are independent, its characteristic triplet is such that $Q=\{q_{i,j}\}$ is a diagonal matrix, i.e.
$$
q_{ii}\geq 0, \quad q_{i,j}=0, \qquad i\neq j, \quad i,j=1,2,...,d,
$$
and the support of $\nu(\dd y)$ is contained in the positive half-axes of $\mathbb{R}^d$, see \cite{Sato} p.67. 
On the $i^{th}$ positive  half-axis 
\begin{gather}\label{rozlozenie ny}
\nu(\dd y)=\nu_i(dy_i),\qquad y=(y_1,y_2,...,y_d),
\end{gather}
for $i=1,2,...,d$.
 The $i^{th}$ coordinate of $Z$ is of infinite variation if and only if its Laplace exponent \eqref{Laplace Zi} is such that $q_{ii}>0$ or
\begin{gather}\label{rerere}
\int_{0}^{1}y_i\nu_i(\dd y_i)=+\infty,
\end{gather}
see \cite[Lemma 2.12]{Kyprianou}. It follows from \eqref{mult. CIR condition} that
$$
\frac{1}{2}\langle QG(x), G(x)\rangle=\frac{1}{2}\sum_{j=1}^{d}q_{jj}G_j^2(x)=cx,
$$
so if $q_{ii}>0$ then $G_i(0)=0$. If it is not the case, using \eqref{rozlozenie ny} and  \eqref{nu G0 finite variation} we see that the integral
\begin{align*}
\int_{(0,+\infty)}v\nu_{G(0)}(\dd v) & =\int_{\mathbb{R}^d_{+}}\langle G(0),y\rangle \nu(\dd y) \\ 
&  =\sum_{j=1}^{d}\int_{(0,+\infty)}G_j(0)y_j \ \nu_j(\dd y_j) =\sum_{j=1}^{d}
G_j(0) \int_{(0,+\infty)}y_j \ \nu_j(\dd y_j),
\end{align*}
is finite, so if  \eqref{rerere} holds then $G_i(0)=0$.\hfill $\square$

 \vskip2ex

\noindent
{\bf Proof of Theorem \ref{TwNiez}:} By assumption \eqref{ass1} and Proposition \ref{rem o G(0)=0} or by assumption \eqref{ass2}  we have $G(0)=0$, so it follows from Remark \ref{rem warunki w eksp. Laplacea} that
\begin{equation} \label{eq:dwa_zero}
J_{Z^{G(x)}}(b) = J_{1}(bG_{1}(x))+J_2(bG_2(x))+...+J_d(bG_d(x))=x\tilde{J}_{\mu}(b), \quad b,x\geq 0,
\end{equation} 
where $\tilde{J}_{\mu}(b) = c b^2 + {J}_{\mu}(b)$, $c\ge 0$ and ${J}_{\mu}(b)$ is given by \eqref{def Jmu i Jnu0}.
This yields
\begin{equation}
\frac{J_{1}\rbr{b\cdot G_{1}(x)}}{J_{1}\rbr{G_{1}(x)}}\cdot\frac{J_{1}\rbr{G_{1}(x)}}{x}+\ldots+\frac{J_{d}\rbr{b\cdot G_{d}(x)}}{J_{d}\rbr{G_{d}(x)}}\cdot\frac{J_{d}\rbr{G_{d}(x)}}{x}=\tilde{J}_{\mu}(b),\label{eq:dwa}
\end{equation}
where in the case $G_i(x) = 0$ we set $\frac{J_{i}\rbr{b\cdot G_{i}(x)}}{J_{i}\rbr{G_{i}(x)}}\cdot\frac{J_{i}\rbr{ G_{i}(x)}}{x} = 0$. Without loss of generality we may assume that $J_{1}$, $J_{2}$,$\ldots$,$J_{d}$ are non-zero (thus positive for positive arguments).  
By assumption, $J_{i}$, $i=1,2,\ldots, d$ vary regularly at $0$ with some indices $\alpha_{i}$, $i=1,2,\ldots,d$, 
so for $b>0$
\begin{equation} \label{eq:trzy}
\lim_{y\ra0+}\frac{J_{i}\rbr{b\cdot y}}{J_{i}(y)}=b^{\alpha_{i}}.
\end{equation}
Assume that 
\[
\alpha_{1}=\ldots=\alpha_{i\rbr{1}}>\alpha_{i\rbr{1}+1}=\ldots=\alpha_{i\rbr{2}}>\ldots\ldots>\alpha_{i\rbr{g-1}+1}=\ldots=\alpha_{i\rbr{g}}=\alpha_{d},
\]
where $i(g) = d$. Let us denote $i_0 = 0$ and
\begin{equation} \label{limits}
\eta_{k}(x) :=\frac{J_{i\rbr{k-1}+1}\rbr{G_{i\rbr{k-1}+1}(x)}+\ldots+J_{i\rbr{k}}\rbr{G_{i\rbr{k}}(x)}}{x}, \quad k=1,2,\ldots, g.
\end{equation}
We can rewrite equation \eqref{eq:dwa} in the form
\begin{equation}
\sum_{k=1}^{g} \rbr{\sum_{i=i\rbr{k-1}+1}^{i\rbr{k}} \frac{J_{i}\rbr{b\cdot G_{i}(x)}}{J_{i}\rbr{G_{i}(x)}}\cdot\frac{J_{i}\rbr{ G_{i}(x)}}{x}}=\tilde{J}_{\mu}(b).\label{eq:trzyy}
\end{equation}
By passing to the limit as $x\ra0+$, from \eqref{eq:trzy} and \eqref{eq:trzyy} we get 
\begin{align}
b^{\alpha_{i\rbr{1}}}\rbr{ \lim_{x \ra 0+} \eta_{1}(x)} +\ldots+   b^{\alpha_{i\rbr{g}}} \rbr{\lim_{x \ra 0+} \eta_{g}(x)}  =\tilde{J}_{\mu}(b), \label{eq:trzyyy}
\end{align} 
thus
\begin{gather}\label{J mu tilde sum power}
\tilde{J}_{\mu}(b)=\sum_{k=1}^g \eta_{k}b^{\alpha_{i\rbr{k}}},
\end{gather}
provided that the limits $\eta_{k} := \lim_{x \ra 0+} \eta_{k}(x)$, $k=1,2,\ldots, g$, exist. 
Thus it remains to prove that for  $k=1,2,\ldots, g$ the limits $\lim_{x \ra 0+} \eta_{k}(x)$ indeed exist {and that $\alpha_{i(k)} \in (1,2]$.}

First we will prove that $\lim_{x \ra 0+} \eta_{g}(x)$ exists.
Assume, by contrary, that this is not true, so 
\begin{equation} 
\limsup_{x \ra 0+} \eta_{g}(x) - \liminf_{x \ra 0+} \eta_{g}(x) \ge \delta >0.
\label{sequencess}
\end{equation}
It follows from \eqref{eq:dwa_zero} that
\begin{equation} \label{ogr}
\frac{J_1(G_1(x))+J_2(G_2(x))+...+J_d(G_d(x))}{x} = \sum_{k=1}^g \eta_k(x)=\tilde{J}_{\mu}(1).
\end{equation}
Let now $b_0 \in (0,1)$ be small enough so that 
\begin{equation} \label{oszacowanie}
\tilde{J}_{\mu}(1) b_0^{\alpha_{i\rbr{g-1}} - \alpha_{i(g)}} < \frac{\delta}{6}.
\end{equation}
Let us set in \eqref{eq:trzyy} $b=b_0$ and then divide both sides of \eqref{eq:trzyy} by $b_0^{\alpha_{i(g)}}$. 
It follows from \eqref{ogr} that each term $\frac{J_{i}\rbr{ G_{i}(x)}}{x}$, $i=1,2,\ldots,d$, is bounded by 
 $\tilde{J}_{\mu}(1)$. From this and \eqref{eq:trzy} for  $x>0$ sufficiently close to $0$ we have 
\[
\eta_g(x) - \frac{\delta}{6} \le \frac{1}{b_0^{\alpha_{i(g)}}} \rbr{ \sum_{i=i\rbr{g-1}+1}^{i\rbr{g}} \frac{{J}_{i}\rbr{b_0\cdot G_{i}(x)}}{J_{i}\rbr{G_{i}(x)}}\cdot\frac{J_{i}\rbr{ G_{i}(x)}}{x}} \le \eta_g(x) + \frac{\delta}{6}
\]
and 
\begin{align*}
\frac{1}{b_0^{\alpha_{i(g)}}}  \sum_{k=1}^{g-1}  \rbr{ \sum_{i=i\rbr{k-1}+1}^{i\rbr{k}} \frac{J_{i}\rbr{b_0\cdot G_{i}(x)}}{J_{i}\rbr{G_{i}(x)}}\cdot\frac{J_{i}\rbr{ G_{i}(x)}}{x}}  \le \sum_{k=1}^{g-1}  2 b_0^{\alpha_{i(k)} - \alpha_{i(g)}} \eta_k(x) \\
\le 2 b_0^{\alpha_{i(g-1)} - \alpha_{i(g)}} \tilde{J}_{\mu}(1)
\end{align*}
thus from \eqref{eq:trzyy}, two last estimates and \eqref{oszacowanie}
\[
\eta_g(x) - \frac{\delta}{6} \le \frac{\tilde{J}_{\mu}(b_0)}{b_0^{\alpha_{i(g)}}} \le \eta_g(x) + \frac{\delta}{6}+ 2\tilde{J}_{\mu}(1) b_0^{\alpha_{i(g-1)} - \alpha_{i(g)}} < \eta_g(x) + \frac{\delta}{2}.
\]
But this contradicts \eqref{sequencess}
since we must have 
\[
 \limsup_{x \ra 0+} \eta_g(x) \le \frac{\tilde{J}_{\mu}(b_0)}{b_0^{\alpha_{i(g)}}} + \frac{\delta}{6}, \quad \liminf_{x \ra 0+} \eta_g(x) \ge \frac{\tilde{J}_{\mu}(b_0)}{b_0^{\alpha_{i(g)}}} - \frac{\delta}{2}.
\]

Having proved the existence of the limits $\lim_{x \ra 0+} \eta_{g}(x)$, ..., $\lim_{x \ra 0+} \eta_{g-m+1}(x)$ we can proceed similarly to prove the existence of the limit $\lim_{x \ra 0+} \eta_{g-m}(x)$. 
Assume that $\lim_{x \ra 0+} \eta_{g-m}(x)$ does not exist, so \begin{equation} 
\limsup_{x \ra 0+} \eta_{g-m}(x) - \liminf_{x \ra 0+} \eta_{g-m}(x) \ge \delta >0.
\label{sequencess1}
\end{equation}
Let $b_0 \in (0,1)$ be small enough so that 
\begin{equation} \label{oszacowanie1}
\tilde{J}_{\mu}(1) b_0^{\alpha_{i\rbr{g-m-1}} - \alpha_{i(g-m)}} < \frac{\delta}{8}.
\end{equation}
Let us set in \eqref{eq:trzyy} $b=b_0$ and then divide both sides of \eqref{eq:trzyy} by $b_0^{\alpha_{i(g-m)}}$. 
For  $x>0$ sufficiently close to $0$ we have 
\[
\eta_{g-m}(x) - \frac{\delta}{8} \le \frac{1}{b_0^{\alpha_{i(g-m)}}} \sum_{i=i\rbr{g-m-1}+1}^{i\rbr{g-m}} \frac{{J}_{i}\rbr{b_0\cdot G_{i}(x)}}{J_{i}\rbr{G_{i}(x)}}\cdot\frac{J_{i}\rbr{ G_{i}(x)}}{x} \le \eta_{g-m}(x) + \frac{\delta}{8},
\]
\begin{align*}
\frac{1}{b_0^{\alpha_{i(g-m)}}}  \sum_{k=1}^{g-m-1}  \rbr{ \sum_{i=i\rbr{k-1}+1}^{i\rbr{k}} \frac{J_{i}\rbr{b_0\cdot G_{i}(x)}}{J_{i}\rbr{G_{i}(x)}}\cdot\frac{J_{i}\rbr{ G_{i}(x)}}{x}}  \le \sum_{k=1}^{g-m-1}  2 b_0^{\alpha_{i(k)} - \alpha_{i(g-m)}} \eta_k(x) \\
\le 2 b_0^{\alpha_{i(g-m-1)} - \alpha_{i(g-m)}} \tilde{J}_{\mu}(1)
\end{align*}
and
\begin{align*}
\sum_{k=g-m+1}^{g}  \frac{b_0^{\alpha_{i(k)}} \eta_k}{b_0^{\alpha_{i(g-m)}}} - \frac{\delta}{8} & \le \frac{1}{b_0^{\alpha_{i(g-m)}}} \sum_{k=g-m+1}^{g}  \sum_{i=i\rbr{k-1}+1}^{i\rbr{k}} \frac{{J}_{i}\rbr{b_0\cdot G_{i}(x)}}{J_{i}\rbr{G_{i}(x)}}\cdot\frac{J_{i}\rbr{ G_{i}(x)}}{x} \\
& \le \sum_{k=g-m+1}^{g}  \frac{b_0^{\alpha_{i(k)}} \eta_k}{b_0^{\alpha_{i(g-m)}}}  + \frac{\delta}{8}
\end{align*}
thus from \eqref{eq:trzyy}, last three estimates and \eqref{oszacowanie1}
\begin{align*}
\eta_{g-m}(x) - \frac{\delta}{4} & \le \frac{J_{\mu}(b_0)}{b_0^{\alpha_{i(g-m)}}} - \sum_{k=g-m+1}^{g}  \frac{b_0^{\alpha_{i(k)}} \eta_k}{b_0^{\alpha_{i(g-m)}}} \\
 & \le \eta_{g-m}(x) + \frac{\delta}{4}+ 2\tilde{J}_{\mu}(1) b_0^{\alpha_{i(g-1)} - \alpha_{i(g)}} <\eta_{g-m}(x) +  \frac{\delta}{2}.
\end{align*}
But this contradicts \eqref{sequencess1}.

Now we are left with the proof that for $k=1,2,\ldots,g$, $\alpha_{i(k)} \in (1,2]$. Since the Laplace exponent of $Z_i$ is given by 
\eqref{Laplace Zi}, by Proposition \ref{bounds_alpha} we necessarily have that $J_i$ varies regularly with index $\alpha_i \in [1,2], i=1,2,...,d$. Thus it remains to prove that $\alpha_i >1, i=1,2,...,d$. If it was not true we would have $\alpha_{i(g)}=1$ in \eqref{J mu tilde sum power} and $\eta_g>0$. Then
$$
\lim_{b\ra 0+} \tilde{J}_{\mu}(b)/b = \lim_{b\ra 0+} J_{\mu}(b)/b =\eta_{g}>0,
$$ 
but, again, by Proposition \ref{bounds_alpha} it is not possible. 
\hfill $\square$

\vskip2ex
\noindent
{\bf Proof of Corollary \ref{cor o postaci mu} :} From Remark \ref{rem warunki w eksp. Laplacea} and Theorem \ref{TwNiez} we know that 
\[
J_{Z^{G(x)}}(b) =  x  c b^2 + x{J}_{\mu}(b) = x \sum_{k=1}^g\eta_{k}b^{\alpha_{{k}}}, 
\]
where $1\le g \le d$, $\eta_{k}> 0$, $\alpha_k\in(1,2]$, $\alpha_k\neq\alpha_j$, $k,j=1,2,\ldots,g$, $c\ge 0$.
Without loss of generality we may assume that $2 \ge \alpha_1 > \alpha_2 > \ldots > \alpha_g >1$. Thus, since the Laplace exponent is nonnegative, $ x{J}_{\mu}(b)$ is of the form
\begin{gather}\label{pierwsza postac}
x{J}_{\mu}(b) = x\sum_{k=1}^g\eta_{k}b^{\alpha_{{k}}}, \qquad \text {if} \ c=0,
\end{gather}
or
\begin{gather}\label{druga postac}
x{J}_{\mu}(b) = x\sbr{(\eta_1-c)b^2+ \sum_{k=2}^g\eta_{k}b^{\alpha_{{k}}}}, \qquad \text {if} \ 0<c\leq\eta_1 \ \text {and} \ \alpha_1=2.
\end{gather}
In the case \eqref{pierwsza postac} we need to show that $\alpha_1<2$. If it was not true, we would have
\[
\lim_{b \ra +\ns} \frac{{J}_{\mu}(b)}{b^2} = \eta_1 >0,
\]
but this contradicts Proposition \ref{bounds_alpha}. In the same way we prove that $\eta_1=c$ in \eqref{druga postac}.
This proves the required representation \eqref{J mu postaaaccccc}.
\hfill $\square$

\subsection{Moments and tails of short rates from the class $\mathbb{A}_g(a,b;\alpha_1,...,\alpha_g; \eta_1,...,\eta_g)$}\label{tail_fatness}

In this section we will prove that  moments of order $p$ of short rates $R(t)$, $t \in (0, +\ns)$, from the class $\mathbb{A}_g(a,b;\alpha_1,...,\alpha_g; \eta_1,...,\eta_g)$ are finite for $p \in \rbr{0,\alpha_g}$ but $R(t)$, $t \in (0, +\ns)$, have fat tails if $\alpha_g \in (1,2)$ and $R(0) = r_0>0$ or $b>0$ in the sense that  then for any $\varepsilon >0$
\[
\E R(t)^{\alpha_g +\varepsilon} = +\ns.
\] 
We will also give some estimates of $\E \rbr{R(t)}^p$ for $p \in \rbr{1,\alpha_g}$. 

Motivated by the form of canonical representations \eqref{rownanie sklajane}  we focus now on the equation
\begin{gather}\label{rownanie numeryka}
dR(t)=(aR(t-)+b)\dd t+\sum_{i=1}^{g}d_i^{1/\alpha_i}R(t-)^{1/\alpha_i}\dd Z^{\alpha_i}(t), \quad R(0)=r_0, \ t>0,
\end{gather}
where $a\in\mathbb{R}, b\geq 0, d_i>0$ and $Z^{\alpha_i}$ is a canonical $\alpha_i$-stable martingale with 
$2\geq \alpha_1>\alpha_2>...>\alpha_g>1$ and $g\geq 1$. By Proposition \ref{prop canonical representation}, \eqref{rownanie numeryka} is the canonical representation of the class $\mathbb{A}_{g}(a,b;\alpha_1,...,\alpha_g;\eta_1,...,\eta_g)$ where 
\begin{gather}\label{ety w kanonicznym}
 d_i = \eta_i/c_{\alpha_i}, \quad i=1,2,\ldots, g,
\end{gather}
and $c_{\alpha_i}$ is given by \eqref{c alpha}. 

The generator of $R$, that is the generator for the solution of \eqref{rownanie numeryka},
takes the form
\begin{align} \label{gennn}
\mathcal{A}f(x)=cx f^{\prime\prime}(x)&+\cbr{a x+b}f^{\prime}(x)+x\int_{(0,+\infty)} \cbr{f(x+v)-f(x)-f^{\prime}(x)v}{\mu}(\dd v),
\end{align}
where
\begin{gather}\label{postac muuu}
{\mu}(\dd v):= \frac{d_l}{v^{1+\alpha_l}} \dd v +...+\frac{d_g}{v^{1+\alpha_g}} \dd v, \quad v>0.
\end{gather}
Recall, if $\alpha_1=2$, then $c=d_1/2$ and $l=2$. Otherwise $c=0$ and $l=1$.

\subsubsection{Moments of the rates $R(t)$, $t \in (0, +\ns)$}

The very first observation we can make  is that the expectation of the solution of \eqref{rownanie numeryka} is equal to $${\cal E}(t) := \E R(t) = \begin{cases} e^{at} {r_0}  + \frac{b}{a} \rbr{e^{at}-1} &\text{ if } a \neq 0,\\ r_0 + bt & \text{ if } a = 0,\end{cases}$$
which readily follows from the fact that it satisfies 
$\dd {\mathcal E}(t) = \rbr{a {\mathcal E}(t) + b} \dd t$, ${\mathcal E}(0) = r_0$. It is also in place to notice that using the product rule for stochastic differentials one checks that $e^{-at}\rbr{R(t) - \E R(t)}$, $t\ge0$, is a martingale. Below we construct another martingale to prove the following result giving a bound for the $p$-th moment of $R(t)$, $p\in(1,\alpha_g)$. 
\begin{prop}\label{moments_r}  If $R(t)$, $t \in (0, +\ns)$, is from the class $\mathbb{A}_g(a,b;\alpha_1,...,\alpha_g; \eta_1,...,\eta_g)$ and $p \in \rbr{0,\alpha_g}$ then $\E \rbr{R(t)}^p < +\ns$. Moreover, for $p \in \rbr{1,\alpha_g}$,
\[
\E \rbr{R(t)}^p \le e^{apt} I_p(t),
\]
where $I_p(t)$ satisfies the following ordinary differential equation
\begin{equation} \label{eqi}
\frac{\dd I_p(t)}{\dd t} =e^{-at}(c(p-1)+b)p\rbr{I_{p}(t)}^{\rbr{p-1}/p} + \sum_{i=l}^{g} e^{a\rbr{1-\alpha_{i}}t} h_{\alpha_i}\rbr{{I}_{p}(t)}^{\rbr{p+1-\alpha_{i}}/p}, \, {I}_{p}(0) = r_0^p,
\end{equation}
with $h_{\alpha_i} = \frac{p(p-1)}{\Gamma(2-p)} \Gamma\rbr{\alpha_i-p} \eta_i$, $i=l,\ldots,g$, where  $l=2$, $c=\eta_1$ in the case when $\alpha_1=2$ while  $l=1$,  $c=0$ in the case when $\alpha_1<2$.
\end{prop}
\noindent 
{\bf Proof:}  
First, we will prove that for any $p\in(1,\alpha_{g})$,
\[
\E(R(t))^{p}<+\ns.
\]
To prove this, let us fix some $B\ge r_{0}$ and $\Delta>0$, and
consider a process $\bar{R}(t)$ which satisfies the equation 
\begin{gather}\label{rownanie numerykabar}
d\bar{R}(t)=(a\bar{R}(t-)+b)\dd t+\sum_{i=1}^{g}d_i^{1/\alpha_i}\bar{R}(t-)^{1/\alpha_i}\dd \bar{Z}^{\alpha_i}(t), \quad \bar{R}(0)=r_0, \ t>0,
\end{gather}
where $\bar{Z}^{\alpha_i}$ are L\'evy martingales with the L\'evy measure
\begin{gather}\label{postac muuubar}
{\bar{\mu}}(\dd v):= \frac{d_l}{v^{1+\alpha_l}}  {\bf 1}_{\cbr{v \le \Delta}} \dd v +...+\frac{d_g}{v^{1+\alpha_g}} {\bf 1}_{\cbr{v \le \Delta}} \dd v, \quad v>0,
\end{gather}
with the same $d_i$, $i=l, \ldots, g$, as in \eqref{ety w kanonicznym}.
Next, let $\tilde{R}(t)$ be the process $\bar{R}(t)$ stopped at
the moment when it reaches the level $B$ or higher, that is 
\[
\tilde{R}(t)=\bar{R}(t\wedge\tau^{B}),\quad t\ge0,
\]
where 
\[
\tau^{B}=\inf\cbr{t\ge0:\bar{R}(t)\ge B}\ge0.
\]
The process $\tilde{R}$ is bounded, thus $\E\left(\tilde{R}(t)\right)^{p}<+\ns$.
Let $f_{p}(r):=r^{p}$ and let $\tilde{{\cal A}}$ be the generator
of $\tilde{R}(t)$. Since $f_{p}$ is convex and increasing, ${\cal \tilde{{\cal A}}}f_{p}(x)$,
$x>0$, is bounded as below
\begin{align}
{\cal \tilde{{\cal A}}}f_{p}(x) & \le cxf_{p}''(x)+(|a|x+b)f_{p}'(x)+x\int_{(0,+\ns)}\cbr{f_{p}(x+v)-f_{p}(x)-f'_{p}(x)v}\mu\rbr{\dd v},\label{eq:boundcalA}
\end{align}
where $\mu\rbr{\dd v}$ is given by \eqref{postac muuu} (recall the generator \eqref{gennn}
of $R$).

For $\alpha\in(p,2)$ we easily calculate 
\begin{align*}
 & \int_{(0,+\ns)}\cbr{f_{p}(x+v)-f_{p}(x)-f'_{p}(x)v}\frac{1}{v^{1+\alpha}}\dd v\\
 & =\int_{(0,+\ns)}\cbr{(x+v)^{p}-x^{p}-px^{p-1}v}\frac{1}{v^{1+\alpha}}\dd v\\
 & =\frac{x^{p}}{x^{\alpha}}\int_{(0,+\ns)}\cbr{(1+u)^{p}-1-pu}\frac{1}{u^{1+\alpha}}\dd u=c_{\alpha,p}x^{p-a},
\end{align*}
where 
\begin{equation}
c_{\alpha,p}:=\frac{\Gamma(2-\alpha)}{\alpha(\alpha-1)}\frac{p(p-1)}{\Gamma(2-p)}\Gamma(\alpha-p).\label{eq:calphap}
\end{equation}
For $x\ge0$ let us define 
\begin{align*}
\tilde{{\cal H}}_{p}(x): & =(c(p-1)+b)px^{p-1}+\sum_{i=l}^{g}c_{\alpha_{i},p}d_{i}x^{p+1-\alpha_{i}}+|a|px^{p}=\\
 & =\sum_{l-1\le i\le g}h_{\alpha_{i}}x^{p+1-\alpha_{i}}+h_{1}x^{p}
\end{align*}
($1<\alpha_{g}<\alpha_{g-1}<\ldots<\alpha_{l}<\alpha_{l-1}=2$, $h_1$, $h_{\alpha_{i}}$ are
defined by the last relation). By (\ref{eq:boundcalA}),
\begin{equation}
{\cal \tilde{{\cal A}}}f_{p}(x)\le\tilde{{\cal H}}_{p}(x).\label{eq:boundcalAH}
\end{equation}
The difference
\[
f_{p}(\tilde{R}(t))-\int_{0}^{t}{\cal \tilde{{\cal A}}}f_{p}\rbr{\tilde{R}(s)}\dd s
\]
is a martingale. For $0<q<\alpha_{g}$, $t\ge0$, we define
\[
\tilde{\cal E}_{q}(t):=\E f_{q}(\tilde{R}(t))=\E(\tilde{R}(t))^{q}.
\]
By Jensen's inequality, for $q<p$
\begin{equation}
\tilde{\cal E}_{q}(t)\le\rbr{\tilde{\cal E}_{p}(t)}^{q/p}.\label{eq:Jensen}
\end{equation}
By (\ref{eq:boundcalAH}), $\tilde{\cal E}_{p}(t)$ satisfies
\begin{align*}
\tilde{\cal E}_{p}(t) & =\E f_{p}(R(t))=\E\int_{0}^{t}{\cal \tilde{{\cal A}}}f_{p}\rbr{\tilde{R}(s)}\dd s=\int_{0}^{t}\E{\cal \tilde{{\cal A}}}f_{p}\rbr{\tilde{R}(s)}\dd s\\
 & \le\int_{0}^{t}\E\tilde{{\cal H}}_{p}(\tilde{R}(s))\dd s=h_{1}\int_{0}^{t}\tilde{\cal E}_{p}(s)\dd s+\sum_{l-1\le i\le g}h_{\alpha_{i}}\int_{0}^{t}\tilde{\cal E}_{p+1-\alpha_{i}}(s)\dd s,
\end{align*}
or, in differential notation,
\begin{equation}
\dd{{\cal E}_{p}(t)}\le h_{1}{\cal E}_{p}(t)\dd t+\sum_{l-1\le i\le g}h_{\alpha_{i}}{\cal E}_{p+1-\alpha_{i}}(t)\dd t.\label{eq:diff_qe}
\end{equation}
Now, using (\ref{eq:Jensen}), we obtain 
\begin{align*}
\dd{\tilde{\cal E}_{p}(t)} & \le h_{1}\tilde{\cal E}_{p}(t)\dd{t+}\sum_{l-1\le i\le g}h_{\alpha_{i}}\tilde{\cal E}_{p+1-\alpha_{i}}(t)\dd t\\
 & \le h_{1}\tilde{\cal E}_{p}(t)\dd t+\sum_{l-1\le i\le g}h_{\alpha_{i}}\rbr{\tilde{\cal E}_{p}(t)}^{\rbr{p+1-\alpha_{i}}/p}\dd t.
\end{align*}
Denoting $h=\sum_{l-1\le i\le g}h_{\alpha_{i}}$ and using the inequality
$y^{\rbr{p+1-\alpha_{i}}/p}\le1+y$ valid for any $y\ge0$, we finally
get the estimate
\[
\dd{\tilde{\cal E}_{p}(t)}\le(h_{1}+h)\rbr{1+\tilde{\cal E}_{p}(t)}\dd t
\]
which yields, that $\tilde{\cal E}_{p}(t)$ is no greater than the solution
of the differential equation
\[
\dd{\tilde E_p(t)} = (h_{1}+h)\rbr{1+\tilde E_p(t)}\dd t,\quad \tilde E_p(0)=r_{0}^{p},
\]
which is equal $e^{(h_{1}+h)t}(r_{0}^p+1)-1$. Thus 
\begin{equation}
\E(\tilde{R}(t))^{p}=\tilde{\cal E}_{p}(t)\le e^{(h_{1}+h)t}(r_{0}^p+1)-1.\label{eq:estimate}
\end{equation}
Let us notice that the estimate (\ref{eq:estimate}) does not depend
on $B$ and $\Delta$. Passing with $B$ and $\Delta$ to $+\ns$
we obtain that $\tilde{R}(t)$ tends almost surely to $R(t)$, thus
$\E(R(t))^{p}<+\ns.$

Now, knowing that $\E(R(t))^{p}<+\ns$ we may reason in a similar
way as before to obtain more precise estimate for $\E(R(t))^{p}$.
Denoting now ${\cal E}_{q}(t):=\E f_{q}(R(t))$ and reasoning in a
similar way as before we obtain the inequality
\begin{equation} \label{bound_growth_of_ep}
\dd{{\cal E}_{p}(t)}\le ap{\cal E}_{p}(t)\dd t+(c(p-1)+b)p\rbr{{\cal E}_{p}(t)}^{\rbr{p-1}/p}\dd t  +\sum_{i=l}^{g}c_{\alpha_{i},p}d_{i}\rbr{{\cal E}_{p}(t)}^{\rbr{p+1-\alpha_{i}}/p} \dd t,
\end{equation}
where $c_{\alpha_{i},p}$ are defined by (\ref{eq:calphap}). Define
\[
{\cal I}_{p}(t):=e^{-apt}{\cal E}_{p}(t),\quad t\ge0
\]
We have 
\begin{align*}
\dd{{\cal I}_{p}(t)} & =-ape^{-apt}{\cal E}_{p}(t)\dd t+e^{-apt}\dd{{\cal E}_{p}(t)}\\
 & \le e^{-apt}(c(p-1)+b)p\rbr{{\cal E}_{p}(t)}^{\rbr{p-1}/p}\dd t + e^{-apt}\sum_{i=l}^{g}c_{\alpha_{i},p}d_{i}\rbr{{\cal E}_{p}(t)}^{\rbr{p+1-\alpha_{i}}/p}\dd t\\
 & =e^{-apt}\sum_{l-1\le i\le g}h_{\alpha_{i}}\rbr{{\cal E}_{p}(t)}^{\rbr{p+1-\alpha_{i}}/p}\dd t\\
 & =e^{-apt}\sum_{l-1\le i\le g}h_{\alpha_{i}}\rbr{e^{apt}{\cal I}_{p}(t)}^{\rbr{p+1-\alpha_{i}}/p}\dd t\\
 & =\sum_{l-1\le i\le g}e^{a\rbr{1-\alpha_{i}}t}h_{\alpha_{i}}\rbr{{\cal I}_{p}(t)}^{\rbr{p+1-\alpha_{i}}/p}\dd t
\end{align*}
($1<\alpha_{g}<\alpha_{g-1}<\ldots<\alpha_{l} < \alpha_{l-1} =2$, $h_{\alpha_{i}}$ are
defined by the last relation). 
This yields, that ${\cal I}_{p}(t)$ is no greater than the solution $I_p(t)$
of the differential equation \eqref{eqi} and
\begin{equation*}
\E({R}(t))^{p}={\cal E}_{p}(t) = e^{apt} {\cal I}_{p}(t) \le e^{apt} { I}_{p}(t).\end{equation*}
\hfill $\square$

\begin{rem}\label{moments_r1}
Let $p \in \rbr{1, \alpha_g}$. Using the notation from the formulation of Proposition \ref{moments_r} and denoting $h=(c(p-1)+b)p+\sum_{l\le i\le g}h_{\alpha_{i}}$, $\beta=a\rbr{1-\alpha_{g}}$,
$\gamma=\rbr{p+1-\alpha_{g}}/p$ we see that if $a \ge 0$ then  ${I}_{p}(t)$ is
no greater than the solution of the differential equation
\[
\frac{\dd{I}(t)}{\dd t}=he^{\beta t}\rbr{I(t)}^{\gamma},\quad I(0)=r_{0}^{p}\vee1,
\]
which is equal
\[
I(t)=\begin{cases}
\rbr{r_{0}^{p(1-\gamma)}\vee1-(1-\gamma)\frac{h}{\beta}\rbr{1-e^{\beta t}}}^{\frac{1}{1-\gamma}} & \text{ if }\beta<0,\\
\rbr{(1-\gamma)ht+r_{0}^{p(1-\gamma)}\vee1}^{\frac{1}{1-\gamma}} & \text{ if }\beta=0.
\end{cases}
\]
This gives that $\E \rbr{R(t)}^p$ grows, as $t \ra +\ns$, no faster than $\text{const}.e^{apt}$ when $a >0$ and no faster than $\text{const}.t^{p/\rbr{\alpha_g - 1}}$ when $a = 0$. 
 
To analyze the situation for $a<0$ let us notice that the function $E_p(t) = e^{apt} I_p(t)$ satisfies the equation  
\begin{equation} \label{pthmoment}
\dd {{E}_{p}(t)} = ap{ E}_{p}(t)\dd t+(c(p-1)+b)p\rbr{{ E}_{p}(t)}^{\rbr{p-1}/p}\dd t  +\sum_{i=l}^{g} h_{\alpha_{i}} \rbr{{E}_{p}(t)}^{\rbr{p+1-\alpha_{i}}/p} \dd t,
\end{equation}
Let $e_p$ be the unique positive solution of the equation
\begin{equation*}
ap \cdot {e}_{p}+(c(p-1)+b)p \cdot {{ e}}^{\rbr{p-1}/p}_{p}  +\sum_{i=l}^{g} h_{\alpha_{i}} {{e}}^{\rbr{p+1-\alpha_{i}}/p}_{p} = 0,
\end{equation*}
If ${ E}_{p}(t) < {e}_{p}$ then $\dd { E}_{p}(t) > 0$ and if ${ E}_{p}(t) > {e}_{p}$ then $\dd { E}_{p}(t) < 0$. From this it follows that $\lim_{t \ra +\ns} { E}_{p}(t) = e_p$ and we obtain $$\limsup_{t \ra +\ns} \E({R}(t))^{p} \le  \lim_{t \ra +\ns} { E}_{p}(t) = e_p$$.
\end{rem}

In what follows we will use the concept of \emph{regularly varying} random vectors introduced in \cite{HultLind07}. For reader's convenience let us recall the definition of such vectors. $\R^g$-valued vector $X$ is regularly varying if there exists a sequence $\rbr{a_n}$ of positive reals such that $a_n \ra +\ns$ and a nonzero Radon measure $\nu$ on the Borel $\sigma$-field ${\cal B} \rbr{ \bar{ \R }^g_0}$ of Borel sets of $\bar{ \R }^g_0 := \rbr{[-\ns, +\ns]^g}\setminus\cbr{0}$ such that 
\begin{equation} \label{regularly varying vectors}
\nu\rbr{[-\ns, +\ns]^g \setminus \R^g} = 0 \text{ and } n \cdot \P\rbr{a_n^{-1}X \in \cdot} \ra^v \nu(\cdot),
\end{equation}
where $ \ra^v $ denotes the vague convergence on ${\cal B} \rbr{ \bar{ \R }^g_0}$. It can be shown that \eqref{regularly varying vectors} implies that there exists some $\alpha >0$ such that for all $u>0$ and $A \in {\cal B} \rbr{ \bar{ \R }^g_0}$ such that $0 \notin\bar{A} $, $\nu \rbr{uA} = u^{-\alpha}\nu(A)$. This is denoted by $X \in RV_{\alpha}\rbr{ \rbr{a_n}, \nu, {\cal B} \rbr{ \bar{ \R }^g_0} }$.

\begin{prop}\label{tails_r}
The rates $R(t)$, $t \in (0, +\ns)$, from the class $\mathbb{A}_g(a,b;\alpha_1,...,\alpha_g; \eta_1,...,\eta_g)$ such that $\alpha_g \in (1,2)$ and $R(0) = r_0>0$ or $b > 0$ have infinite moments of order $\alpha_g +\varepsilon$ for any $\varepsilon >0$. 
\end{prop}
\noindent 
{\bf Proof:}  
Let ${Z(t)}=\rbr{{Z}^{\alpha_1}(t), {Z}^{\alpha_2}(t),...,{Z}^{\alpha_g}(t)}$, $t \ge 0$, be vector of canonical stable martingales with indices $\alpha_k, k=1,2,...,g$, respectively. Using \eqref{tails_Zzz} it is easy to notice that $Z(1) \in RV_{\alpha_g}\rbr{ \rbr{a_n}, \nu, {\cal B} \rbr{ \bar{ \R }^g_0} }$ with $a_n = n^{1/\alpha_g}$ and the $\alpha_g$-stable measure $\nu$ concentrated on the $g$th half-axis:
\[
\nu\rbr{\dd z_1, \dd z_2,  \ldots, \dd z_g} = \delta_0(\dd z_1) \ldots \delta_0(\dd z_{g-1}) \frac{1}{\alpha_g}\frac{1}{z_g^{1+\alpha_g}} {\bf 1}_{\cbr{z_g>0}} \dd z_g,
\]
where $\delta_0$ denotes Dirac's delta measure on $\R$ concentrated at $0$.
By \eqref{rownanie numeryka} $R(t)$ has the same distribution as the stochastic integral $r_0 + \rbr{Y\cdot \tilde{Z}}(t)$ with the predictable c\`adl\`ag integrand $$Y(t) = \rbr{aR(t-)+b, d_1^{1/\alpha_1}R(t-)^{1/\alpha_1}, \ldots, d_g^{1/\alpha_g}R(t-)^{1/\alpha_g}}$$ and the integrator $\tilde{Z}(t) = \rbr{t, {Z}^{\alpha_1}, {Z}^{\alpha_2},...,{Z}^{\alpha_g}}$. Assume that there exists some $\varepsilon >0$ such that
\begin{equation} \label{wrong}
\E R(t)^{\alpha_g +\varepsilon} < +\ns.
\end{equation} 
Since $e^{-at}\rbr{R(t) - \E R(t)}$, $t\ge0$, is a martingale,
by the Doob maximal $L^p$ inequality applied to this martingale we obtain
\[
\E \rbr{ \sup_{0\le s\le t}R(s)}^{\alpha_g +\varepsilon} < +\ns.
\]
This and the form of the integrand $Y$ means that we can apply  \cite[Theorem 3.4]{HultLind07} and obtain that $$R(t) \in RV_{\alpha_g}\rbr{ \rbr{a_n}, \nu^*, {\cal B} \rbr{ \bar{ \R }_0} },$$ where the measure $\nu^*$ does not vanish. But this yields that for any $\varepsilon >0$
\[
\E R(t)^{\alpha_g +\varepsilon} = +\ns.
\] 
which is a contradiction with \eqref{wrong}. \hfill $\square$

\begin{rem}\label{tails_r1}
We conjecture that when the assumptions of Proposition \ref{tails_r} are satisfied then in fact $\E R(t)^{\alpha_g} = +\ns$ as it is the case for stable CIR models ( $\mathbb{A}_1(a,b;\alpha_1; \eta_1)$ models with $\alpha_1 \in (1,2)$ in our notation), see \cite[Proposition 3.1]{LiMa}. However, for the brevity of the proof we decided to restrict to considering the moments strictly greater than $\alpha_g$. 
\end{rem}

\subsubsection{Limit distributions of the rates $R(t)$ as $t \ra +\ns$ and their tails}

General results on the limit distributions of CBI processes are proven in \cite{Li}, see also \cite{MiKe}, \cite{Pinsky} and \cite[Proposition 3.7]{JiaoMaScotti} (however, in our opinion the statement of \cite[Proposition 3.7]{JiaoMaScotti} is not true for all '$\alpha$-CIR integral type processes' defined in \cite{JiaoMaScotti} since even an $\alpha$-CIR process may not possess the limit distribution).

To state the condition on the existence of the limit distributions of the rates $R(t)$, $t \in (0, +\ns)$, from the class $\mathbb{A}_g(a,b;\alpha_1,...,\alpha_g; \eta_1,...,\eta_g)$ we shall define two functions $\mathcal{R}, \mathcal{F}$ (\emph{branching mechanism} and \emph{immigration mechanism}, respectively). $\mathcal{R}$ and $\mathcal{F}$ depend on the generator of $R$ and are defined as
\begin{align}\label{wzor na F}\nonumber
\mathcal{R}(\lambda)&:=-c\lambda^2+\Big[a+\int_{(1,+\infty)}(1-v){\mu}(\dd v)\Big]\lambda+\int_{0}^{+\infty}(1-e^{-\lambda v}-\lambda(1\wedge v)){\mu}(\dd v),\\[1ex]
\mathcal{F}(\lambda)&:=b\lambda.
\end{align} 
From \eqref{postac muuu} we obtain
 \begin{align}\label{wzor na R}\nonumber
\mathcal{R}(\lambda)&=-c\lambda^2+a\lambda-\int_{0}^{+\infty}(e^{-\lambda v}-1+\lambda v){\mu}(\dd v)\\[1ex]
& =-c\lambda^2+a\lambda-\sum_{i=l}^g\eta_k\lambda^{\alpha_k}=a\lambda-\sum_{i=1}^g\eta_k\lambda^{\alpha_k} .
\end{align}

By \cite[Theorem 2.6]{MiKe} the following statements are equivalent:
\begin{itemize}
\item $R(t)$, $t\ge0$, converges (as $t \ra +\ns$) in distribution to some random variable $R_{\ns}$ with the distribution $\mathcal L$;
\item $R(t)$, $t \ge 0$, has the unique invariant distribution $\mathcal L$;
\item it holds that $a = \mathcal{R}'(0) \le 0$ and
\[
-\int_0^u \frac{\mathcal{F}(\lambda)}{\mathcal{R}(\lambda)} \dd \lambda <+\ns
\]
for some $u>0$.
\end{itemize}
Moreover, the limit distribution $\mathcal L$, in the case it exists, is infinitely divisible and its Laplace transform reads 
\[
\E \exp\rbr{-u R_{\ns}} = \exp \rbr{\int_0^u \frac{\mathcal{F}(\lambda)}{\mathcal{R}(\lambda)} \dd \lambda}, \quad u\ge 0.
\]
From these statements we obtain the following.

\begin{prop}\label{lim_r}
The rates $R(t)$, $t \in (0, +\ns)$, from the class $\mathbb{A}_g(a,b;\alpha_1,...,\alpha_g; \eta_1,...,\eta_g)$ converge (as $t \ra +\ns$) in distribution to some random variable $R_{\ns}$ iff one of the following holds: (i) $b=0$ and $a\le0$ or (ii) $b>0$, $a<0$ or (iii) $b>0$, $a=0$ and $\alpha_g <2$. In the case (i) $\P \rbr{R_{\ns} = 0} = 1$, in the case (ii) $\E R_{\ns} = -b/a$ and in the case (iii) the tail of $R_{\ns}$ has the asymptotics
\[
\P\rbr{R_{\ns} > r} \sim \frac{b}{\eta_g\rbr{2-\alpha_g} \Gamma\rbr{\alpha_g - 1}} \frac{1}{r^{2-\alpha_g}} \text{ as } r \ra +\ns.
\]
\end{prop}
{\bf Proof:}  

(i) In this case  the ratio $-{\mathcal{F}(\lambda)}/{\mathcal{R}(\lambda)}$ reduces to $0$ and the statement follows. 

(ii) In this case the ratio $-{\mathcal{F}(\lambda)}/{\mathcal{R}(\lambda)}$ reduces to 
\[
-\frac{\mathcal{F}(\lambda)}{\mathcal{R}(\lambda)} = \frac{b}{-a + \sum_{i=1}^g\eta_k\lambda^{\alpha_k-1}} = -\frac{b}{a} + o(1) \text{ as } \lambda \ra 0+
\] 
and we obtain
\begin{align}
\E \exp\rbr{-u R_{\ns}} = \exp \rbr{\int_0^u \frac{\mathcal{F}(\lambda)}{\mathcal{R}(\lambda)} \dd \lambda} =  1+\frac{b}{a}u + o(u) \text{ as } u \ra 0+.
\end{align}
Hence, $$
\E R_{\ns} = \lim_{u \ra 0+} \frac{\E \exp\rbr{-u R_{\ns}}-1}{-u} = -\frac{b}{a}.
$$

(iii) In this case the ratio $-{\mathcal{F}(\lambda)}/{\mathcal{R}(\lambda)}$ reduces to 
\[
-\frac{\mathcal{F}(\lambda)}{\mathcal{R}(\lambda)} = \frac{b}{\sum_{i=1}^g\eta_k\lambda^{\alpha_k-1}} = \frac{b}{\eta_g \lambda^{\alpha_g - 1}} + o\rbr{ \frac{1}{\lambda^{\alpha_g - 1}}} \text{ as } \lambda \ra 0+
\] 
and we obtain
\begin{align}
& 1-\exp \rbr{\int_0^u \frac{\mathcal{F}(\lambda)}{\mathcal{R}(\lambda)} \dd \lambda} = 1-\exp \rbr{-\int_0^u \frac{b}{\eta_g}  \lambda^{1-\alpha_g}+ o\rbr{{\lambda^{1-\alpha_g }}} \dd \lambda} \\
& = 1-\rbr{1-\frac{b}{\eta_g\rbr{2-\alpha_g}}u^{2-\alpha_g} + o\rbr{u^{2-\alpha_g}}} = \frac{b}{\eta_g\rbr{2-\alpha_g}} u^{2-\alpha_g} + o\rbr{u^{2-\alpha_g}} \text{ as } u \ra 0+.
\end{align}
Hence, by the Tauberian theorem \cite[Corollary 8.1.7]{Bingham}, 
\[
\P\rbr{R_{\ns} > r} \sim \frac{b}{\eta_g\rbr{2-\alpha_g} \Gamma\rbr{\alpha_g - 1}} \frac{1}{r^{2-\alpha_g}} \text{ as } r\ra +\ns.
\]

\hfill $\square$

\subsection{Generating equations on a plane}\label{section Generalized CIR equations on a plane}

In this section we characterize all equations \eqref{rownanie 2}, with $d=2$, which generate affine models by 
a direct description of the classes $\mathbb{A}_1(a,b; \alpha_1;\eta_1)$ and $\mathbb{A}_2(a,b;\alpha_1,\alpha_2; \eta_1,\eta_2)$.
Our analysis requires an additional regularity assumption that the components of $G$ are strictly positive outside zero and 
\begin{gather}\label{iloraz G regularny}
	\frac{G_2(\cdot)}{G_1(\cdot)}\in C^1(0,+\infty).
\end{gather}
Then $\mathbb{A}_1(a,b; \alpha_1;\eta_1)$ consists of the following equations\\
$$
\bullet \quad \dd R(t)=(aR(t)+b)\dd t+c_0 R(t)^{1/\alpha_1}\Big(G_1 \dd Z_1(t)+G_2 \dd Z_2(t)\Big), 
$$
where $c_0=({\eta_1}/{c_{\alpha_1}})^{1/{\alpha_1}}$, $G_1,G_2$ are positive constants and $G_1 Z_1(t)+G_2 Z_2(t)$ is an $\alpha_1$-stable process,
$$
\bullet \quad \dd R(t)=(aR(t)+b)\dd t+G_1(R(t-))\dd Z_1(t)+\left(\frac{\eta_1 R(t-)-c_1 G_1^{\alpha_1}(R(t-))}{c_2}\right)^{1/\alpha_1}\dd Z_2(t), 
$$
where $c_1,c_2>0$, $G_1(\cdot)$ is any function such that  
$$
G_1(x)> 0, \quad \frac{\eta_1 x-c_1 G_1^{\alpha_1}(x)}{c_2}>0, \qquad x>0,
$$
and $Z_1, Z_2$ are stable processes with index $\alpha_1$.\\
\noindent
The class $\mathbb{A}_2(a,b;\alpha_1,\alpha_2; \eta_1,\eta_2)$ is a singleton.

The classification above follows directly from the following result.

\begin{tw}\label{tw d=2 independent coord.}
Let $G(x)=(G_1(x), G_2(x))$ be continuous functions such that $G_1(x)>0,G_2(x)>0, x>0$ and \eqref{iloraz G regularny} holds.
Let $Z(t)=(Z_1(t),Z_2(t))$ have independent coordinates of infinite variation with Laplace exponents varying regularly at zero
with indices $\alpha_1,\alpha_2$, respectively, where $2\geq \alpha_1\geq\alpha_2>1$. 
\begin{enumerate}[I)]
\item If $\tilde{J}_{\mu}$  is of the form 
\begin{gather}\label{pierwszy przyp Jmu}
	\tilde{J}_{\mu}(b)=\eta_1 b^{\alpha_1}, \quad b\geq 0,
\end{gather}
with $\eta_1>0, 1<\alpha_1\leq 2$, then $(G,Z)$ is a  generating pair if and only if one of the following two cases holds:
\begin{enumerate}[a)]
\item 
\begin{gather}\label{wspolinionwosc G}
G(x)=c_0 \ x^{1/\alpha_1}\cdot\left(
\begin{array}{ccc}
 G_1\\
 G_2, 
\end{array}
\right), \quad x\geq 0,
\end{gather}
where $c_0=(\frac{\eta_1}{c_{\alpha_1}})^{\frac{1}{\alpha_1}}, G_1>0, G_2>0$ and the process
$$
G_1 Z_1(t)+G_2 Z_2(t), \quad t\geq 0,
$$
is $\alpha_1$-stable.
\item $G(x)$ is such that
\begin{gather}\label{Ib}
c_1 G^{\alpha_1}_1(x)+c_2 G^{\alpha_1}_2(x)=\eta_1 x, \quad x\geq 0,
\end{gather}
with some constants $c_1,c_2>0$, and $Z_1,Z_2$ are $\alpha_1$-stable processes.
\end{enumerate}

\item If $\tilde{J}_{\mu}$  is of the form 
\begin{gather}\label{drugi przyp Jmu}
	\tilde{J}_{\mu}(b)=\eta_1 b^{\alpha_1}+\eta_2 b^{\alpha_2},\quad b\geq 0,
\end{gather}

with $\eta_1,\eta_2>0, 2\geq \alpha_1>\alpha_2>1$ then $(G,Z)$ is a  generating pair if and only if
\begin{gather}\label{postac g w drugim prrzypadku}
G_1(x)=\left(\frac{\eta_1}{c_1}  x\right)^{1 / \alpha_1}, \quad G_2(x)=\left(\frac{\eta_2}{d_2}  x\right)^{1 / \alpha_2}, \quad x\geq 0,
\end{gather}
with some $c_1,d_2>0$ and $Z_1$ is $\alpha_1$-stable, $Z_2$ is $\alpha_2$-stable.
\end{enumerate}
\end{tw}

For the proof of Theorem \ref{tw d=2 independent coord.} we refer to Sect. \ref{proof_d2} of the Appendix. 

\subsection{An example in 3D}\label{section Example in higher dimensions}

In Section \ref{section Generalized CIR equations on a plane} we proved that in the case $d=2$ the set $\mathbb{A}_2(a,b;\alpha_1,\alpha_2;\eta_1,\eta_2)$ is a singleton. Here we show that  
this property breaks down when $d=3$. In the example below we construct a family of generating pairs 
$(G,Z)$ such that 
\begin{gather}
	J_{Z^{G(x)}}(b)=x\left(\eta_1 b^{\alpha_1}+\eta_2 b^{\alpha_2}\right),\quad b\geq 0,
\end{gather}
with $\eta_1,\eta_2>0, 2\geq \alpha_1>\alpha_2>1$ and such that the related generating equations differ from the canonical representation of $\mathbb{A}_2(a,b;\alpha_1,\alpha_2;\eta_1,\eta_2)$.

\begin{ex} 
Let us consider a process $Z(t)=(Z_1(t),Z_2(t),Z_3(t))$ with independent coordinates such that  $Z_1$ is $\alpha_1$-stable, $Z_2$ is $\alpha_2$-stable,  $Z_3$ is a sum of an $\alpha_1$-  and $\alpha_2$-stable processes. Then
$$
J_1(b)=\gamma_1 b^{\alpha_1}, \quad J_2(b)=\gamma_2 b^{\alpha_2}, \quad J_3(b)=\gamma_3b^{\alpha_1}+\tilde{\gamma}_3 b^{\alpha_2}, \quad b\geq 0,
$$
where $\gamma_1>0,\gamma_2>0,\gamma_3>0,\tilde{\gamma}_3>0$. We are looking for non-negative functions $G_1, G_2,G_3$ solving the equation
\begin{gather}\label{krokodyl}
J_1(bG_1(x))+J_2(bG_2(x))+J_3(bG_3(x))=x \left(\eta_1 b^{\alpha_1}+\eta_2 b^{\alpha_2}\right), \quad x,b\geq 0.
\end{gather}
It follows from \eqref{krokodyl} that
$$
\gamma_1b^{\alpha_1}(G_1(x))^{\alpha_1}+\gamma_2b^{\alpha_2}(G_2(x))^{\alpha_2}
+\gamma_3b^{\alpha_1}(G_3(x))^{\alpha_1}+\tilde{\gamma}_3 b^{\alpha_2}(G_3(x))^{\alpha_2}=
x\left[\eta_1 b^{\alpha_1}+\eta_2b^{\alpha_2}\right], \quad x,b\geq 0,
$$
and, consequently,
$$
b^{\alpha_1}\left[\gamma_1G_1^{\alpha_1}(x)+\gamma_3 G_3^{\alpha_1}(x)\right]+
b^{\alpha_2}\left[\gamma_2G_2^{\alpha_2}(x)+\tilde{\gamma}_3 G_3^{\alpha_2}(x)\right]=x\left[\eta_1 b^{\alpha_1}+\eta_2b^{\alpha_2}\right], \quad x,b\geq 0.
$$
Thus we obtain the following system of equations 
\begin{gather*}
\gamma_1G_1^{\alpha_1}(x)+\gamma_3 G_3^{\alpha_1}(x)=x\eta_1, \\
\gamma_2G_2^{\alpha_2}(x)+\tilde{\gamma}_3 G_3^{\alpha_2}(x)=x\eta_2,
\end{gather*}
which allows us to determine $G_1$ and $G_2$ in terms of $G_3$, that is
\begin{gather}\label{G1}
G_1(x)=\left(\frac{1}{\gamma_1}\left(x\eta_1-\gamma_3 G_3^{\alpha_1}(x)\right)\right)^{\frac{1}{\alpha_1}}\\\label{G2}
G_2(x)=\left(\frac{1}{\gamma_2}\left(x\eta_2-\tilde{\gamma}_3 G_3^{\alpha_2}(x)\right)\right)^{\frac{1}{\alpha_2}}.
\end{gather}
The positivity of $G_1,G_2,G_3$  means that $G_3$ satisfies
\begin{gather}\label{war nieuj}
0\leq G_3(x)\leq \left(\frac{\eta_1}{\gamma_3}x\right)^{\frac{1}{\alpha_1}}\wedge \left(\frac{\eta_2}{\tilde{\gamma}_3}x\right)^{\frac{1}{\alpha_2}}, \quad x\geq 0.
\end{gather}
It follows that $(G,Z)$ with any $G_3$ satisfying \eqref{war nieuj} and $G_1,G_2$ given by \eqref{G1}, \eqref{G2} constitutes a generating pair.
\end{ex}

\section{Applications}\label{section Applications}

In this section, We investigate the relevance of the equation \eqref{rownanie numeryka} to the description of risk-free market rates. 
First, in Section \ref{section Bond prices in canonical models}, we describe the dependence of the arising bond prices 
$$
P(t,T)=e^{-A(T-t)-B(T-t)R(t)}, 
$$
 on the parameters by describing the dependence
\begin{align*}
A(t,T)&=A(t,T)(a,b,d_1,...,d_g,\alpha_1,...,\alpha_g),\\
B(t,T)&=B(t,T)(a,b,d_1,...,d_g,\alpha_1,...,\alpha_g).
\end{align*}
Then, in Section \ref{section Calibration of canonical models}, we pass to the calibration of the resulting model rates to the rate quotes of the European Central Bank. The source of data we use can be found  at:
\newline
\url{
https://www.ecb.europa.eu/stats/financial_markets_and_interest_rates/euro_area_yield_curves/html/index.en.html}.
\newline\noindent
It covers a wide time range $2004-2024$ embracing the whole spectrum of states of the  European economy.
The resulted variety of the market data  allows us to test and to judge the performance of the model in a reliable way.
In particular, we compare the model generated by \eqref{rownanie numeryka} with a standard CIR model.

\subsection{Bond prices in canonical models}\label{section Bond prices in canonical models}

Let us start with recalling the concept of pricing based on the semigroup
\begin{gather}\label{semigroup}
\mathcal{Q}_tf(x):= \mathbb{E}[e^{-\int_{0}^{t}R(s)ds}f(R(t))\mid R(0)=x], \quad t\geq 0,
\end{gather}
which was developed in \cite{FilipovicATS}. The formula provides the price at time $0$ of the claim $f(R(t))$ paid at time $t$ given $R(0)=x$. By Theorem 5.3 in \cite{FilipovicATS} for $f_{\lambda}(x):=e^{-\lambda x}, \lambda \geq 0$ we know that
\begin{gather}\label{wycena ogolnie}
\mathcal{Q}_t f_\lambda (x)=e^{-\rho(t,\lambda)-\sigma(t,\lambda)x}, \quad x\geq 0,
\end{gather}
where $\sigma(\cdot,\cdot)$ satisfies the equation
$$
\frac{\partial\sigma}{\partial t}(t,\lambda)=1+ \mathcal{R}(\sigma(t,\lambda)), \quad \sigma(0,\lambda)=\lambda,
$$
and $\rho(\cdot,\cdot)$ is given by
$$
\rho(t,\lambda)=\int_{0}^{t}\mathcal{F}(\sigma(s,\lambda))ds,
$$
where $\mathcal{R}$ and $\mathcal{F}$ are defined in \eqref{wzor na R} and \eqref{wzor na F}.

Application of the pricing procedure above for $f_\lambda$ with $\lambda=0$  
allows us to obtain from \eqref{wycena ogolnie} the prices of zero-coupon bonds. Using the closed form formula 
\eqref{wzor na R} leads to the following result.
\begin{tw}\label{tw o cenach obligacji}
The zero-coupon bond prices in the affine model generated by \eqref{rownanie numeryka} are equal
\begin{gather}\label{affine prices tw}
	P(t,T)=e^{-A(T-t)-B(T-t)R(t)},
\end{gather}
where $B$ and $A$ are such that 
	\begin{align}\label{B indep. coord.}
	B^\prime(v)&=1+aB(v)-\sum_{i=1}^{g}\eta_iB^{\alpha_i}(v),\quad B(0)=0,\\[1ex]\label{A indep. coord.}
	A^\prime(v)&=b B(v), \quad A(0)=0,
	\end{align}
with $\eta_i$, $i=1, \ldots, g$, given by \eqref{ety w kanonicznym}.
\end{tw}

In the case when $g=1$ and $\alpha_1=2$ equation \eqref{B indep. coord.} becomes a Riccati equation and its explicit solution
provides bond prices for the classical CIR equation. In the opposite case \eqref{B indep. coord.} can be solved 
by numerical methods which exploit the tractable form of the function $\mathcal{R}$ given by \eqref{wzor na R}. Note that $\mathcal{R}$ is continuous, $\mathcal{R}(0)=1$ and  $\lim_{\lambda\rightarrow +\infty}\mathcal{R}(\lambda)=-\infty$. Thus 
$\lambda_0:=\inf\{\lambda>0:1+ \mathcal{R}(\lambda)=0\}$
is a positive number and 
\begin{gather}\label{zachowanie R w lambda 0}
1+ \mathcal{R}(\lambda_0)=0, \quad \mathcal{R}^{\prime}(\lambda_0)<0. 
\end{gather}
The function
\begin{gather}\label{funkcja G}
\mathcal{G}(x):=\int_{0}^{x}\frac{1}{1+ \mathcal{R}(y)}dy, \quad x\in[0,\lambda_0),
\end{gather}
is strictly increasing and its behaviour near $\lambda_0$ can be estimated by
substituting $z=\frac{1}{\lambda_0-y}$ in \eqref{funkcja G} and using the inequality
$$
(\lambda_0-h)^{\alpha}\geq \lambda_0^{\alpha}-\alpha\lambda_0^{\alpha-1}h, \quad h\in(0, \lambda_0),\quad \alpha\in(1,2) .
$$
For the case when $\alpha_1=2$ this yields for $x \in [0, \lambda_0)$
\begin{align}\label{oszacowanie dla G}\nonumber
\mathcal{G}(x)&=\int_{1/\lambda_0}^{1/(\lambda_0-x)}\frac{1}{1+\mathcal{R}(\lambda_0-\frac{1}{z})}\cdot\frac{1}{z^2} \dd z\\[1ex]\nonumber
&=\int_{1/\lambda_0}^{1/(\lambda_0-x)}\frac{ \dd z}{z^2+a\lambda_0 z^2-az-\eta_1(\lambda_0 z-1)^2-\sum_{i=2}^{g}\eta_i z^2(\lambda_0-\frac{1}{z})^{\alpha_i}} \ \\[1ex]\nonumber
&\geq \int_{1/\lambda_0}^{1/(\lambda_0-x)}\frac{ \dd z}{z^2+a\lambda_0 z^2-az-\eta_1(\lambda_0 z-1)^2-\sum_{i=2}^{g}\eta_i z^2(\lambda_0^{\alpha_i}-\alpha_i\lambda_0^{\alpha_i-1}\frac{1}{z})} \ \\[1ex]\nonumber
&=\int_{1/\lambda_0}^{1/(\lambda_0-x)}\frac{\dd z}{z^2(1+a\lambda_0-\eta_1\lambda_0^2-\sum_{i=2}^{g}\eta_i\lambda_0^{\alpha_i})+z(2\eta_1\lambda_0-a+\sum_{i=2}^{g}\alpha_i\eta_i\lambda_0^{\alpha_i-1})-\eta_1} \ \\[1ex]
&=\int_{1/\lambda_0}^{1/(\lambda_0-x)}\frac{\dd z}{(1+\mathcal{R}(\lambda_0))z^2-\mathcal{R}^{\prime}(\lambda_0)z-\eta_1} \ .\end{align}
It follows from \eqref{oszacowanie dla G} and \eqref{zachowanie R w lambda 0} that
$$
\lim_{x\rightarrow \lambda_0^{-}}\mathcal{G}(x)=+\infty,
$$ 
so $\mathcal{G}$ is invertible and $\mathcal{G}^{-1}$ exists on $[0,+\infty)$. Writing \eqref{B indep. coord.} as
$$
B^\prime(v)=1+\mathcal{R}(B(v)), \quad B(0)=0,
$$
we see that 
$$
\frac{d}{dv}\mathcal{G}(B(v))=\frac{1}{1+\mathcal{R}(B(v))}B^\prime(v)=1,
$$
and consequently
$$
\mathcal{G}(B(v))=v, \quad v\geq 0.
$$
Representing $B(\cdot)$ as the inverse of $\mathcal{G}(\cdot)$ enables its numerical computation. 

\subsection{Calibration of canonical models}\label{section Calibration of canonical models}
Our calibration procedure is concerned with the spot rates 
\begin{gather}\label{spot rates empirical}
\hat{y}(T_i), \ i=1,2,...,N,
\end{gather}
of European Central Bank (ECB) which are computed from the zero coupon AAA-rated bonds. The maturity grip $\{T_1,...,T_N\}$ consists of $N=13$ points: $3,6,9$ months and $1,2,3,4,5,10,15,20$, $25$ and $30$ years. Densely chosen small maturities save rapid changes of the market yield curve
$$
T\rightarrow \hat{y}(T),
$$
near zero, while sparsely distributed large maturities do not change the tail shape of the curve. 
The model spot rates are given by 
\begin{gather}\label{spot rates in the model}
y(T_i):=\frac{1}{T_i}\left(\frac{1}{P(0,T_i)}-1\right), \quad i=1,2,...,N,
\end{gather}
where $\{P(0,T_i)\}_i$ denote the bond prices at time $t=0$ generated by the equation \eqref{rownanie numeryka}.
The dependence of the model spot rates on the parameters $(a,b,\alpha_1,...,\alpha_g, d_1,...,d_g)$  is hidden in the function 
$
\mathcal{G}(x)=\int_{0}^{x}{1}/\rbr{1+\mathcal{R}(y)}dy, $
as its inverse enables computing the bond prices via solving the equations (\ref{B indep. coord.}-\ref{A indep. coord.}) for $A(\cdot)$ and $B(\cdot)$.  The calibration aim is to minimize the fitting error measured by a relative distance of the  model spot rates \eqref{spot rates in the model} from the empirical ones \eqref{spot rates empirical}. It is given by the formula
\begin{gather}\label{error spot rates}
Error(a,b,\alpha_1,...,\alpha_g, d_1,...,d_g):=\sum_{i=1}^{N}\frac{(y(T_i)-\hat{y}(T_i))^2}{\hat{y}^2(T_i)}.
\end{gather}
In what follows we compare this error with the error of the CIR model.

\subsubsection{Fitting of the $\alpha$-CIR model to market data}

We start with fitting the $\alpha$-CIR model. Then \eqref{rownanie numeryka} takes the form
\begin{gather}\label{rownanie calibration}
dR(t)=(aR(t)+b)\dd t+d_1^{1/2}R(t-)^{1/2}\dd W(t)+d_2^{1/\alpha}R(t-)^{1/\alpha}\dd Z^{\alpha}(t),
\end{gather}
and the calibration error is minimized with respect to the parameters $(a,b,d_1,d_2,\alpha)$.
The case  $d_2=0$ yields the CIR model. 

Our numerical results of calibration at randomly chosen $15$ dates reveal 
a significant reduction of the calibration error by the $\alpha$-CIR model in most cases as compared to the CIR model.
In over $66\%$ of cases the error reduction exceeds $10\%$. 
In  $46\%$ cases the improvement is greater than $30\%$ and in $33\%$ greater than $50\%$. In $20\%$ cases the error can not be reduced, see Tab. \ref{tabela zbiorcza} for details.

\begin{figure}[htb]
\begin{center}
\begin{tabular}{| c || c | c | c |}   
   \hline
   
\multirow{2}{*}{} & \multirow{2}{*}{Calibration error $\times 100$} &\multirow{2}{*}{Calibration error  $\times 100$}&\multirow{2}{*}{Percentage} \\[1ex]  
Date&in CIR&in $\alpha$-CIR&error reduction\\ \hhline{|====|}
15.07.2008& 0.092 & 0.091	&1.96 \\ \hline
03.06.2009&	15.214	&15.214&	0.00\\ \hline
17.08.2010&	  10.194&	  10.194& 	0.00\\ \hline
06.10.2010&	    2.352&	  0.599&	74.50\\ \hline    
21.10.2011&	    4.712&	  3.289&	30.21\\ \hline
20.03.2012&     12.503&  7.247&	42.04\\ \hline
04.07.2012&	328.212&	  307.701&	6.25\\ \hline
23.09.2013&	  42.196&	  36.007&	14.67\\ \hline
03.12.2014&	177.865&	  56.485&	68.24\\ \hline
03.07.2015&	    1.484&	   1.328&	10.54\\ \hline
10.01.2018&	    0.951&	   0.445&	53.22\\ \hline
21.11.2018&	    0.486&	   0.240&	50.55\\ \hline
21.10.2019&    	13.979&	 13.979&	0.00\\ \hline
08.04.2022&	     24.102&	0.831&	 96.55  \\ \hline
30.08.2022&	   11.499&	8.960& 	22.08\\ \hline  
\end{tabular}
    \captionof{table}{Calibration errors of CIR and $\alpha$-CIR models for the ECB rates at randomly chosen dates.}\label{tabela zbiorcza}
\end{center}
\end{figure}

\begin{figure}[htb]
\begin{center}
\includegraphics[height=2.5in,width=3in,angle=0]{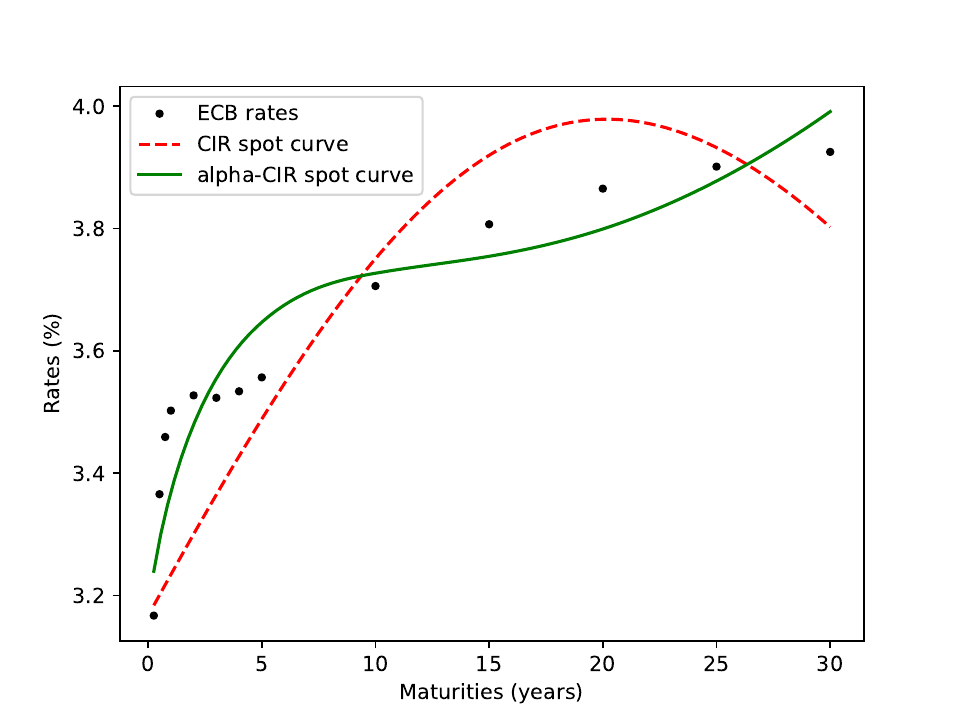}\\
\includegraphics[height=2.5in,width=3in,angle=0]{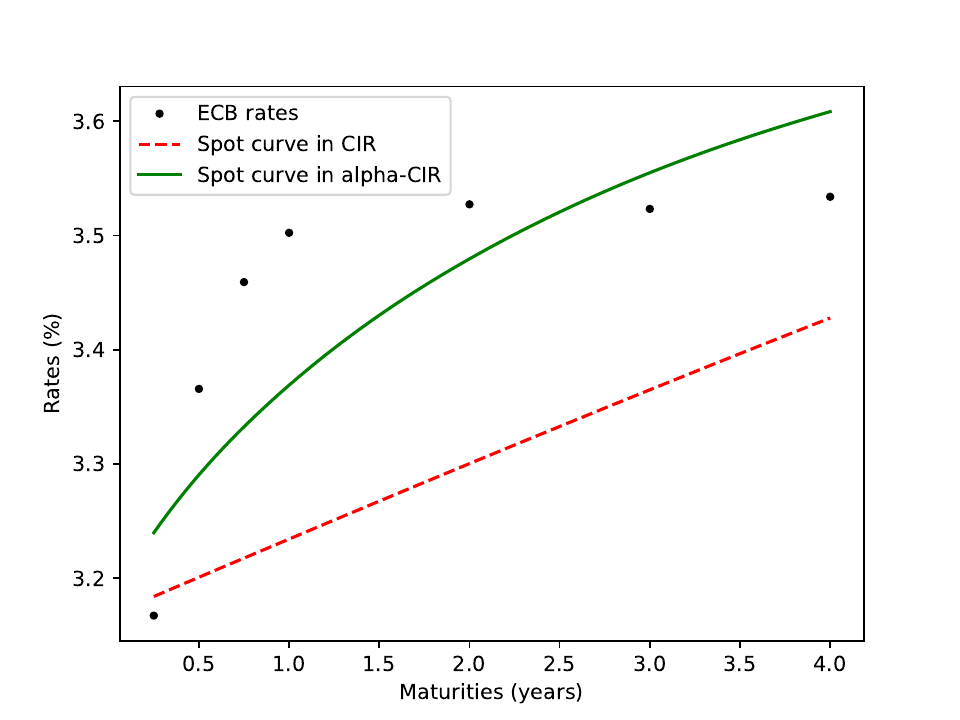}
\includegraphics[height=2.5in,width=3in,angle=0]{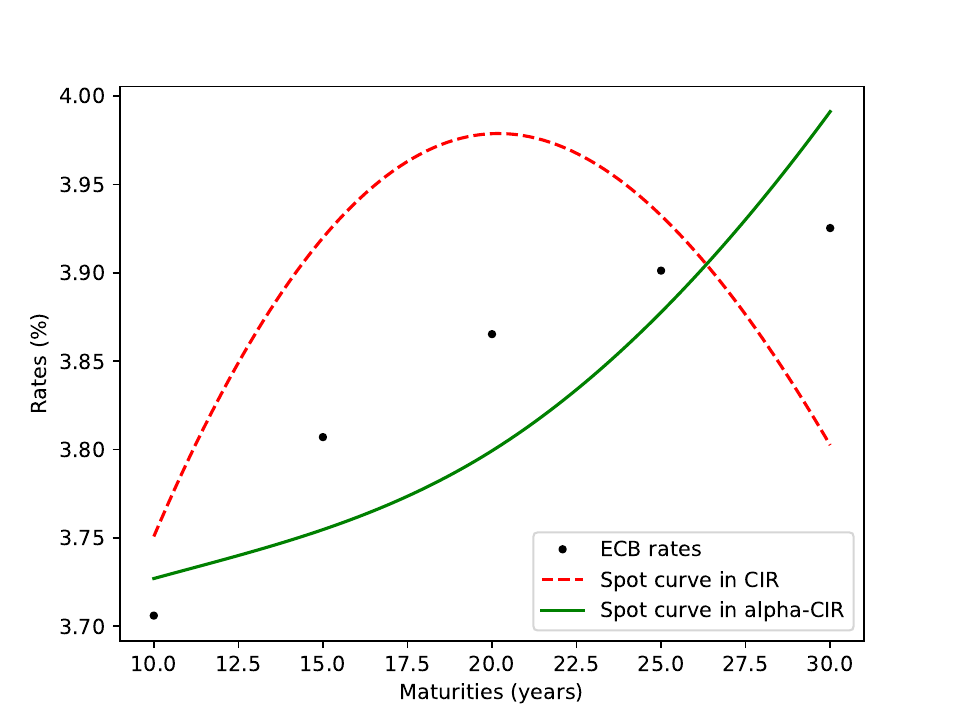}
\caption{Calibration to the ECB rates from 6.10.2010. View for all/small/large maturities.}\label{plot 2010}
\end{center}
\end{figure}

\begin{figure}[htb]
\begin{center}
\includegraphics[height=2.5in,width=3in,angle=0]{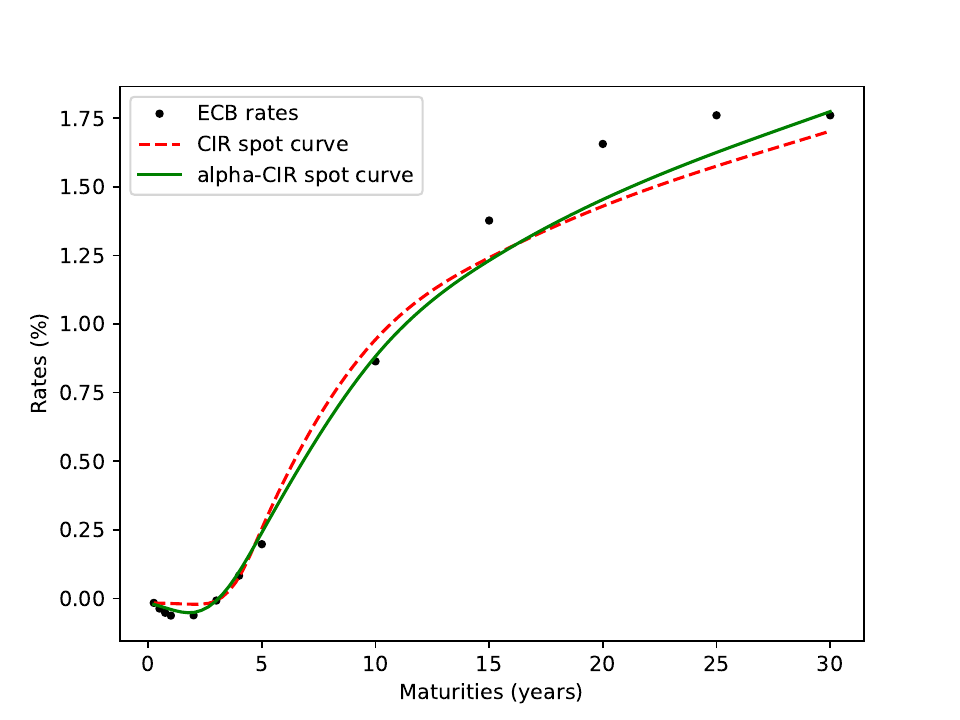}\\
\includegraphics[height=2.5in,width=3in,angle=0]{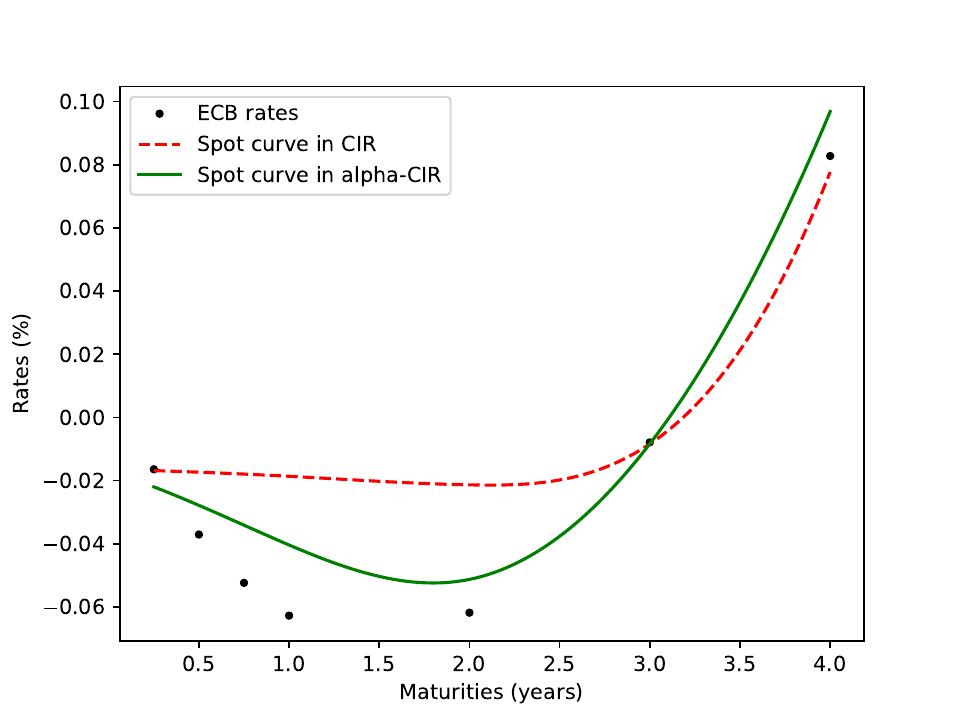}
\includegraphics[height=2.5in,width=3in,angle=0]{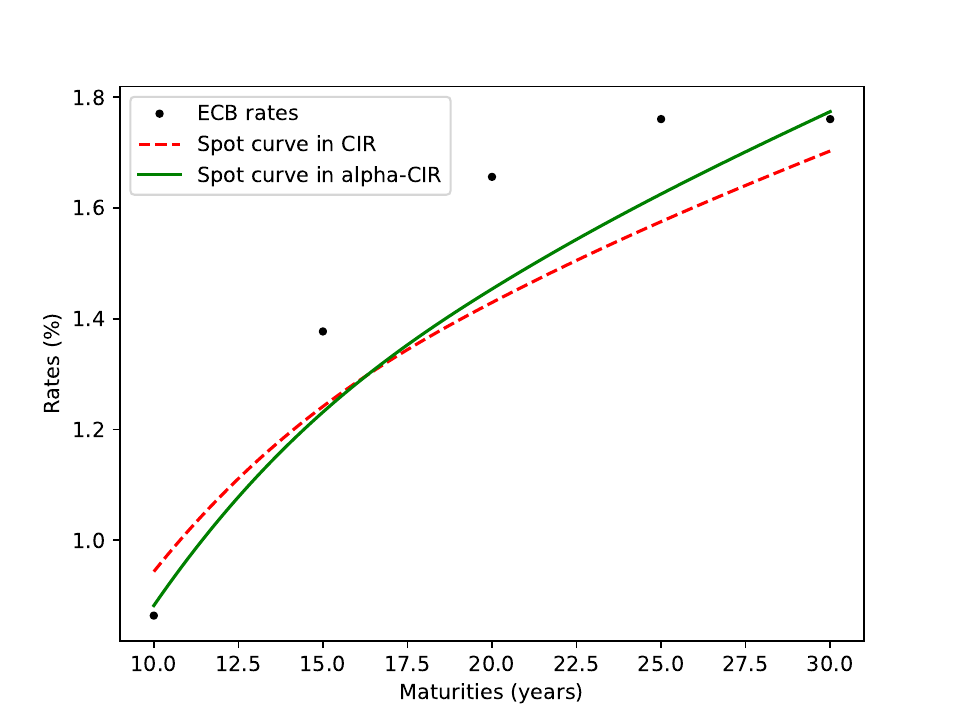}
\caption{Calibration to the ECB rates from 3.12.2014. View for all/small/large maturities.}\label{plot 2014}
\end{center}
\end{figure}

\begin{figure}[htb]
\begin{center}
\includegraphics[height=2.5in,width=3in,angle=0]{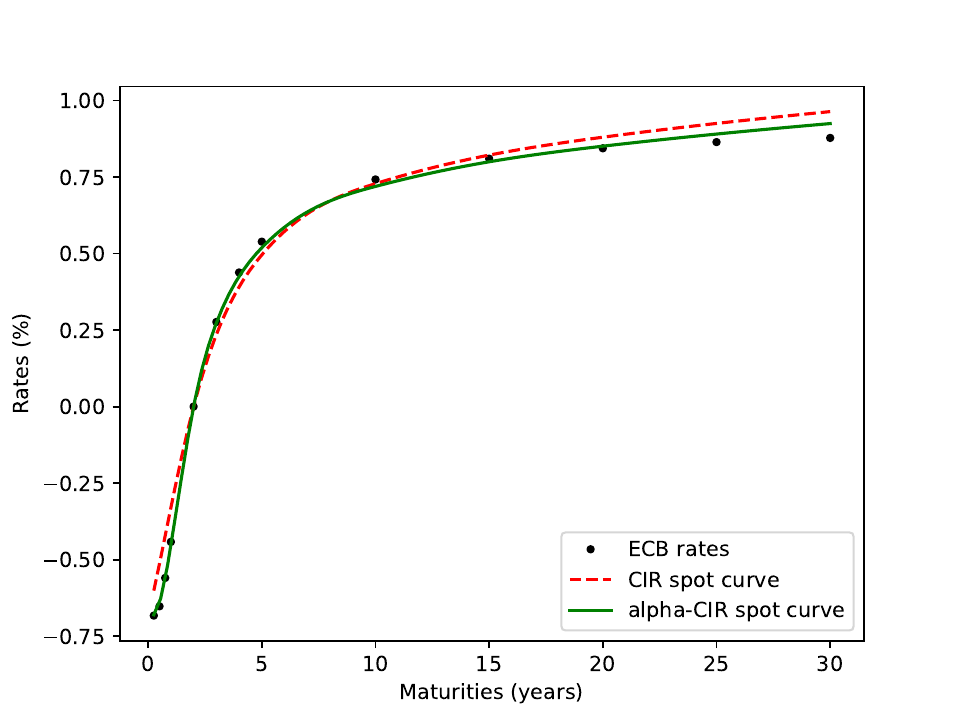}\\
\includegraphics[height=2.5in,width=3in,angle=0]{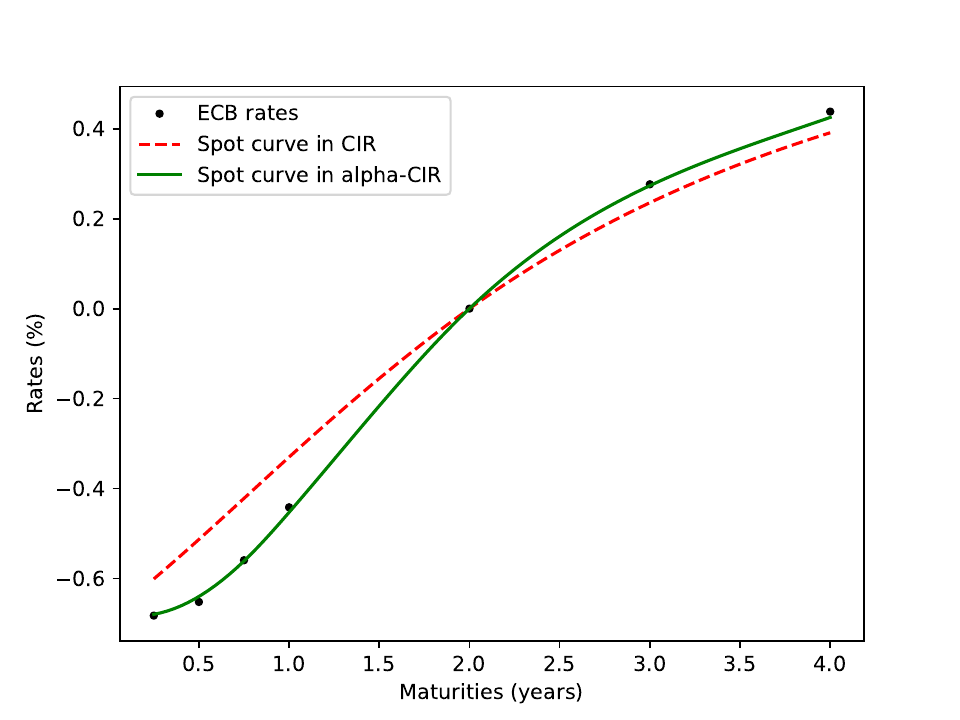}
\includegraphics[height=2.5in,width=3in,angle=0]{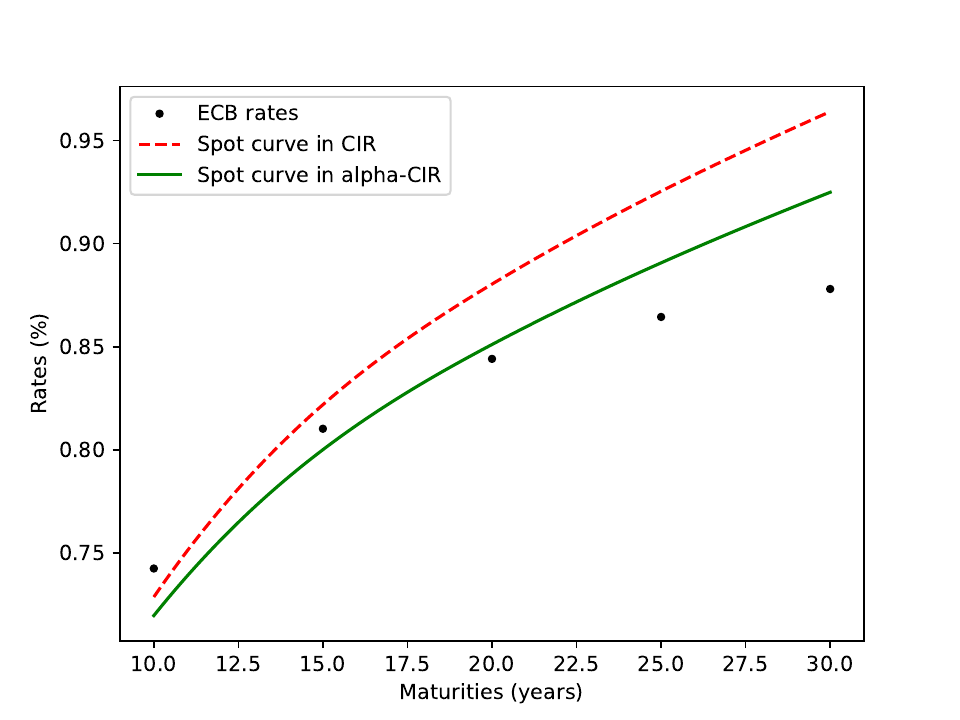}
\caption{Calibration to the ECB rates from 8.04.2022. View for all/small/large maturities.}\label{plot 2022}
\end{center}
\end{figure}

Parameters of models at dates with the best performance are presented in Tab. \ref{tabela parameters} and related plots in 
Fig.\ref{plot 2010}, Fig.\ref{plot 2014} and Fig.\ref{plot 2022}. The yield curves in the $\alpha$-CIR model turn out to be much more flexible than the CIR curves and almost all market quotes are better approached by the $\alpha$-CIR curve.

\begin{figure}[htb]
\begin{center}
\begin{tabular}{|c||c||c|c|c|c|c||c|}
\hline
\multirow{2}{*}{Date} & \multirow{2}{*}{Model} &\multicolumn{5}{c|}{Parameters}& \multirow{2}{*}{Error  $\times 100$} \\  \cline{3-7}
&&$a$&$b$&$d_1$&$d_2$&$\alpha$&\\ \hhline{|========|}
\multirow{2}{*}{06.10.2010} & CIR & 0.011	&0.000&	0.002&& & 2.352\\ \cline{2-8}
&$\alpha$-CIR &   1.384 &0.002 &0.000 &0.074	&1.031&0.599\\ \hline
\multirow{2}{*}{03.12.2014} & CIR & 1.332&	 0.000	& 0.032 && & 177.865\\ \cline{2-8}
 &$\alpha$-CIR &  9.999	& 0.000	& 1.55 x $10^{-6}$ & 0.752	&1.057	& 56.485 \\ \hline
\multirow{2}{*}{08.04.2022} & CIR & -0.066	&0.006	 &0.692 && & 24.102\\ \cline{2-8}
&$\alpha$-CIR &  0.939	&0.005 &1.902&9.12 x $10^{-6}$	&1.999 &0.831\\ \hline
\end{tabular}
 \captionof{table}  {Parameters of CIR and $\alpha$-CIR models for dates with the greatest error reduction.}\label{tabela parameters}
\end{center}
\end{figure}

In our investigation we did not observe a significant error reduction by adding further stable noise components to the equation \eqref{rownanie calibration}. Typical results we
obtained in our implementation look like those in Tab.\ref{table 2010}, Tab.\ref{table 2014} and Tab.\ref{table 2022}, where GCIR($k$) stands for the generalized CIR equation with $k$-dimensional noise.  However, despite negligible calibration improvements, dimensions $k\geq 3$ may be desirable to adjust the heaviness of the tails of $R$. We observed that in the calibration process of GCIR($k$), $k \ge 3$, our algorithm was adding heavier noises than the one $\alpha$-stable noise appearing in the $\alpha$-CIR model.


\newpage
\begin{figure}[h]
\begin{center}
    \begin{tabular}{| c | c | l |}
    \hline    
    Model & Calibration error  $\times 100$ & Stability indices \\ \hline
    CIR & 2.35231805& $\alpha=2$ \\ \hline
    GCIR(1)& 0.59992129&$\alpha_1=1.031$\\ \hline
        GCIR(2)& 0.5907987&$\alpha_1=1.02$, $\alpha_2=1.014$\\ \hline
    GCIR(3)& 0.58845355&$\alpha_1=1.017$, $\alpha_2=1.008$, $\alpha_3=1.007$\\ \hline
   GCIR(4)&0.058819833&$\alpha_1=1.00944$, $\alpha_2=1.00943$, $\alpha_3=1.008$, $\alpha_4= 1.006$\\ \hline
   GCIR(5)&0.058796683 &$\alpha_1=1.01$, $\alpha_2=1.009$, $\alpha_3=1.00898$, $\alpha_4=1.00896$, $\alpha_5=1.0003$\\ \hline
    \end{tabular}
    \captionof{table}{Error reduction - calibration to the ECB rates from 6.10.2010.}\label{table 2010}    
\end{center}
\end{figure}

\newpage
\begin{figure}[h]
\begin{center}
    \begin{tabular}{| c | c | l |}
    \hline    
    Model & Calibration error  $\times 100$ & Stability indices \\ \hline
    CIR & 177.86453440& $\alpha=2$ \\ \hline
         GCIR(1) & 177.86453440& $\alpha=2$ \\ \hline
        GCIR(2)& 56.48507979&$\alpha_1=2$, $\alpha_2=1.057$\\ \hline
    GCIR(3)& 56.44050315&$\alpha_1=2$, $\alpha_2=1.057$, $\alpha_3=1.056$\\ \hline
   GCIR(4)&56.44050315&$\alpha_1=2$, $\alpha_2=1.602$, $\alpha_3=1.05728$, $\alpha_4= 1.05724$\\ \hline
   GCIR(5)&56.44050315 &$\alpha_1=2$, $\alpha_2=1.611$, $\alpha_3=1.589$, $\alpha_4=1.05724$, $\alpha_5=1.05722$\\ \hline
    \end{tabular}
    \captionof{table}{Error reduction - calibration to the ECB rates from 3.12.2014.}\label{table 2014}    
\end{center}
\end{figure}

\begin{figure}[h]
\begin{center}
    \begin{tabular}{| c | c | l |}
    \hline    
    Model & Calibration error  $\times 100$ & Stability indices \\ \hline
    CIR & 24.10280133& $\alpha=2$ \\ \hline
    GCIR(1) & 24.10280133& $\alpha=2$ \\ \hline
    GCIR(2)& 0.83059934&$\alpha_1=2$, $\alpha_2=1.99$\\ \hline
    GCIR(3)& 0.83055904&$\alpha_1=2$, $\alpha_2=1.17$, $\alpha_3=1.14$\\ \hline
   GCIR(4)&0.83050323&$\alpha_1=2$, $\alpha_2=1.35$, $\alpha_3=1.25$, $\alpha_4= 1.21$\\ \hline
   GCIR(5)&0.83049801&$\alpha_1=2$, $\alpha_2=1.53$, $\alpha_3=1.48$, $\alpha_4=1.35$, $\alpha_5=1.23$\\ \hline
    \end{tabular}
    \captionof{table}{Error reduction - calibration to the ECB rates from 8.04.2022.}\label{table 2022}    
\end{center}
\end{figure}

\clearpage
\subsubsection{Remarks on computational methodology}

Our computations were implemented in the Python programming language. The calibration error was minimized with the use of the Nelder-Mead algorithm which turned out to be most effective among all available algorithms for local minimization in the Python library. The computation time of calibration for the $\alpha$-CIR model lied in most cases in the range 100-300 seconds. Calibration of models with a higher number of noise components took typically about 800 seconds but some outliers with 10.000 seconds also appeared.  This stays in a strong contrast  to the CIR model for which the closed form formulas shorten the calibration to the 2 second limit. We suspect that global optimization algorithms would provide even better fit, but they were too slow for the data with more than several maturities.

\clearpage
\newpage

\section{Appendix}

\subsection{Proof of Proposition \ref{prop wstepny}}
\noindent
{\bf Proof:} $(A)$ It was shown in \cite [Theorem 5.3]{FilipovicATS} that the generator of a general positive Markovian short rate generating  an affine model is of the form
\begin{align}\label{generator Filipovica}
\mathcal{A}f(x)=&c x f^{\prime\prime}(x)+(\beta x+\gamma)f^\prime(x)\\[1ex] \nonumber
&+\int_{(0,+\infty)}\Big(f(x+y)-f(x)-f^\prime(x)(1\wedge y)\Big)(m(\dd y)+x\mu(\dd y)), \quad x\geq 0,
\end{align}
for $f\in\mathcal{L}(\Lambda)\cup C_c^2(\mathbb{R}_{+})$, where 
$\mathcal{L}(\Lambda)$ is the linear hull of $\Lambda:=\{f_\lambda:=e^{-\lambda x}, \lambda\in(0,+\infty)\}$
and $C_c^2(\mathbb{R}_{+})$ stands for the set of twice continuously differentiable functions with compact support in $[0,+\infty)$. 
Above $c, \gamma\geq 0$, $\beta\in\mathbb{R}$ and $m(\dd y)$, $\mu(\dd y)$ are nonnegative Borel measures on $(0,+\infty)$ satisfying
\begin{gather}\label{warunki na iary Filipovica}
\int_{(0,+\infty)}(1\wedge y)m(\dd y)+\int_{(0,+\infty)}(1\wedge y^2)\mu(\dd y)<+\infty.
\end{gather}

The generator of the short rate process given by \eqref{rownanie 2} equals
\begin{align*}
 \mathcal{A}_{R}f(x) =& f^\prime(x)F(x)+\frac{1}{2}f^{\prime\prime}(x)\langle QG(x),G(x)\rangle \\
& +\int_{\mathbb{R}^d}\Big(f(x+\langle G(x),y\rangle)-f(x)-f^\prime(x)\langle G(x),y\rangle\Big)\nu(\dd y) \\
 = & f^\prime(x)F(x)+\frac{1}{2}f^{\prime\prime}(x)\langle QG(x),G(x)\rangle \\
&+\int_{\mathbb{R}}\Big(f(x+v)-f(x)-f^\prime(x)v\Big)\nu_{G(x)}(\dd v)
\end{align*}
where  $f$ is a bounded, twice continuously differentiable function. 

By Proposition \ref{prop o skokach z dodatniosci rozwiazania} below, the support of the measure $\nu_{G(x)}$ is contained in $[-x,+\ns)$, thus it follows that  
\begin{align}\label{Generatorr R gen}\nonumber
\mathcal{A}_{R}f(x) = &f^\prime(x)F(x)+\frac{1}{2}f^{\prime\prime}(x)\langle QG(x),G(x)\rangle \\ \nonumber
&+\int_{(0, +\ns)}\Big(f(x+v)-f(x)-f^\prime(x)(1\wedge v) \Big)\nu_{G(x)}(\dd v)\\ \nonumber
&+f^\prime(x)\int_{(0, +\ns)}\Big((1\wedge v)-v\Big)\nu_{G(x)}(\dd v)\\[1ex] \nonumber 
&+\int_{(-\ns, 0)}\Big(f(x+v)-f(x)-f^\prime(x)v \Big)\nu_{G(x)}(\dd v)\\ \nonumber
= & \frac{1}{2}f^{\prime\prime}(x)\langle QG(x),G(x)\rangle + f^\prime(x)\sbr{F(x)+\int_{(1,+\ns)}\Big(1- v\Big)\nu_{G(x)}(\dd v)} \\[1ex] \nonumber
&+\int_{(0, +\ns)}\Big(f(x+v)-f(x)-f^\prime(x)(1\wedge v) \Big)\nu_{G(x)}(\dd v)\\ 
&+\int_{[-x, 0)}\Big(f(x+v)-f(x)-f^\prime(x)v \Big)\nu_{G(x)}(\dd v).  
\end{align}

Comparing \eqref{Generatorr R gen} with \eqref{generator Filipovica} applied to a function $f_{\lambda}$ with $\lambda >0$ such that $f_{\lambda}(x) = e^{-\lambda x}$ for $x \ge 0$, we get

\begin{align}
&c x \lambda^2  - (\beta x+\gamma)  \lambda  \nonumber \\[1ex] \nonumber
&+ \int_{(0,+\infty)}\Big(e^{-\lambda y}-1+\lambda (1\wedge y)\Big)(m(\dd y)+x\mu(\dd y)) \\[1ex] \nonumber
& -  \frac{1}{2} \lambda^2 \langle QG(x),G(x)\rangle + \sbr{F(x)+\int_{(1,+\ns)}\Big(1- v\Big)\nu_{G(x)}(\dd v)} \lambda  
\\[1ex] \nonumber
& - \int_{(0,+\infty)}\Big(e^{-\lambda v}-1+\lambda (1\wedge v)\Big)\nu_{G(x)}(\dd v) \\[1ex] \label{uaua}
& =  \int_{[-x, 0)}\Big(e^{-\lambda v}-1+\lambda  v \Big)\nu_{G(x)}(\dd v) , \quad \lambda >0, x \geq 0.
\end{align}
Comparing the left and the right sides of \eqref{uaua} we see that the left side grows no faster than a quadratic polynomial of $\lambda$ while the right side grows faster that $d e^{\lambda y}$ for some $d, y >0$, unless  the support of the measure $\nu_{G(x)}(\dd v)$ is contained in $[0,+\ns)$. It follows that $\nu_{G(x)}(\dd v)$ is concentrated on $[0,+\infty)$, hence $(a)$ follows, and 

\begin{align}
&c x \lambda^2  - (\beta x+\gamma)  \lambda  \nonumber \\[1ex] \nonumber
& -  \frac{1}{2} \lambda^2 \langle QG(x),G(x)\rangle + \sbr{F(x)+\int_{(1,+\ns)}\Big(1- v\Big)\nu_{G(x)}(\dd v)} \lambda  
\\[1ex] \label{uaua1}
& = \int_{(0,+\infty)}\Big(e^{-\lambda y}-1+\lambda (1\wedge y)\Big)\rbr{\nu_{G(x)}(\dd y) - m(\dd y)- x\mu(\dd y)}, \quad \lambda >0, x \geq 0.
\end{align}
Dividing both sides of the last equality by $\lambda^2$ and using the  estimate  
$$\frac{e^{-\lambda y}-1+\lambda (1\wedge y)}{\lambda^2} \le \rbr{\frac{1}{2} y^2} \wedge \rbr{\frac{e^{-\lambda}-1+\lambda}{\lambda^2}}$$ 
we get that that the left side of \eqref{uaua1} converges to $c x -  \frac{1}{2} \langle QG(x),G(x)\rangle$ as $\lambda\rightarrow +\infty$, while the right side converges to  $0$. This yields \eqref{mult. CIR condition}, i.e. 
\begin{align}\label{W1}
c x=& \frac{1}{2}\langle QG(x),G(x)\rangle,\quad x\geq 0.
\end{align}
Next, fixing $x\geq 0$ and comparing \eqref{Generatorr R gen} with \eqref{generator Filipovica} applied to a function from the domains of both generators and such that $f(x) = f'(x) = f ''(x) =0$ we get 
\[
\int_{(0,+\infty)} f(x+y) (m(\dd y)+x\mu(\dd y)) =  \int_{(0, +\ns)} f(x+v) \nu_{G(x)}(\dd v) 
\]
for any such a function,
which yields 
\begin{gather}\label{rozklad nuG na sume}
\nu_{G(x)}(\dd v)\mid_{(0,+\infty)}=m(\dd v)+x\mu(\dd v), \quad x\geq 0.
\end{gather} 
This implies also
\begin{align}
\label{W2}
\beta x+\gamma=&F(x)+\int_{(1,+\ns)}\Big(1 - v\Big)\nu_{G(x)}(\dd v),\quad x\geq 0.
\end{align}

$(b)$ Setting $x=0$ in \eqref{rozklad nuG na sume} yields 
\begin{gather}\label{nu G0}
\nu_{G(0)}(\dd v)\mid_{(0,+\infty)}=m(\dd v).
\end{gather}
To prove \eqref{nu G0 finite variation}, by \eqref{warunki na iary Filipovica} and \eqref{nu G0}, we need to show that
\begin{gather}\label{cocococ}
\int_{(1,+\infty)}v\nu_{G(0)}(\dd v)<+\infty.
\end{gather}
It is true if $G(0)=0$ and for $G(0)\neq 0$ the following estimate holds
\begin{align*}
\int_{(1,+\infty)}v\nu_{G(0)}(\dd v)&=\int_{\mathbb{R}^d}\langle G(0),y\rangle\mathbf{1}_{[1,+\infty)}(\langle G(0),y\rangle)\nu(\dd y)\\[1ex]
&\leq \mid G(0)\mid\int_{\mathbb{R}^d}\mid y\mid\mathbf{1}_{[1/\mid G(0)\mid,+\infty)}(\mid y\mid)\nu(\dd y),
\end{align*}
and \eqref{cocococ} follows.

$(c)$   \eqref{rozklad nu G(x)} follows from \eqref{rozklad nuG na sume} and \eqref{nu G0}.
To prove \eqref{war calkowe na mu} we use \eqref{rozklad nu G(x)},  \eqref{nu G0 finite variation}
and the following estimate for $x\geq 0$:
\begin{align*}
\int_{0}^{+\infty}(v^2\wedge v)\nu_{G(x)}(\dd v)&=\int_{\mathbb{R}^d}(\mid\langle G(x),y\rangle\mid^2\wedge\langle G(x),y\rangle)\nu(\dd y)\\[1ex]
&\leq \Big(\mid G(x)\mid^2\vee \mid G(x)\mid\Big)\int_{\mathbb{R}^d}(\mid y\mid^2\wedge \mid y\mid)\nu(\dd y)<+\infty.
\end{align*}

$(d)$ It follows from \eqref{W2} and \eqref{rozklad nu G(x)} that 
\begin{align*}
\beta x+\gamma&=F(x)+\int_{(1,+\infty)}(1 - v)\nu_{G(x)}(\dd v)\\[1ex]
&=F(x)+\int_{(1,+\infty)}(1-v)\nu_{G(0)}(\dd v)+x\int_{(1,+\infty)}(1-v)\mu(\dd v), \quad x\geq 0.
\end{align*}
Consequently, \eqref{linear drift}
follows with
$$
a:=\Big(\beta-\int_{(1,+\infty)}(1-v)\mu(\dd v)\Big), \ b:=\Big(\gamma-\int_{(1,+\infty)}(1-v)\nu_{G(0)}(\dd v)\Big),
$$
and $b\geq \int_{(1,+\infty)}(v-1)\nu_{G(0)}(\dd v)$ because $\gamma\geq 0$.

$(B)$ We use \eqref{W2},  \eqref{linear drift} and \eqref{rozklad nuG na sume} to write \eqref{generator Filipovica} in the form
\begin{align*}
\mathcal{A}f(x)=cx f^{\prime\prime}(x)&+\Big[ax +b+\int_{(1,+\infty)}(1 -v)\nu_{G(x)}(\dd v)\Big]f^{\prime}(x)\\[1ex]
&+\int_{(0,+\infty)}[f(x+v)-f(x)-f^{\prime}(x)(1\wedge v)]\nu_{G(x)}(\dd v)\}.
\end{align*}
In view of \eqref{rozklad nuG na sume} and \eqref{nu G0} we see that
\eqref{generator R w tw} is true.

\begin{prop}\label{prop o skokach z dodatniosci rozwiazania}
Let $G:[0,+\infty)\rightarrow \mathbb{R}^d$ be continuous. If the equation \eqref{rownanie 2} has a non-negative strong solution for any initial condition $R(0)=x\geq 0$, then
\begin{gather}\label{ograniczenia skokow}
\forall x\geq 0 \quad \nu{\{y\in\mathbb{R}^d: x+\langle G(x),y\rangle<0\}}=0.
\end{gather}
In particular, the support of the measure $\nu_{G(x)}(\dd v)$ is contained in $[-x,+\infty)$.
\end{prop}
{\bf Proof:} Let us assume to the contrary, that for some $x\geq 0$
$$
\nu{\{y\in\mathbb{R}^d: x+\langle G(x),y\rangle<0\}}>0.
$$
Then there exists $c>0$ such that
$$
\nu{\{y\in\mathbb{R}^d: x+\langle G(x),y\rangle<-c\}}>0.
$$
Let $A\subseteq \{y\in\mathbb{R}^d: x+\langle G(x),y\rangle<-c\}$ be a Borel set separated from zero. By the continuity of $G$ we have that for some $\varepsilon>0$:
\begin{gather}\label{pani w szpileczkach}
\tilde{x}+\langle G(\tilde{x}),y\rangle<-\frac{c}{2}, \quad  \tilde{x}\in[(x-\varepsilon)\vee 0,x+\varepsilon],\quad y\in A.
\end{gather}
Let $Z^2$ be a L\'evy processes with characteristics $(0,0,\nu^2(dy))$, where $\nu^2(dy):=\mathbf{1}_{A}(y)\nu(dy)$ and $Z^1$ be defined by $Z(t)=Z^1(t)+Z^2(t)$. Then $Z^1, Z^2$ are independent and $Z^2$ is a compound Poisson process. Let us consider 
the following equations
\begin{gather*}
dR(t)=F(R(t))dt+\langle G(R(t-)),dZ(t)\rangle, \quad R(0)=x,\\[1ex]
dR^1(t)=F(R^1(t))dt+\langle G(R^1(t-)),dZ^1(t)\rangle, \quad R^1(0)=x.
\end{gather*}
For the exit time $\tau_1$ of $R^1$ from the set $[(x-\varepsilon)\vee 0,x+\varepsilon]$ and the first jump time $\tau_2$ of $Z^2$ we can find
$T>0$ such that $\mathbb{P}(\tau_1>T, \tau_2<T)=\mathbb{P}(\tau_1>T)\mathbb{P}(\tau_2<T)>0$. On the set $\{\tau_1>T, \tau_2<T\}$ we have $R(\tau_2-)=R^1(\tau_2-)$ and therefore 
$$
R(\tau_2)=R^1(\tau_2-)+\langle G(R^1(\tau_2-)),\triangle Z^2(\tau_2)\rangle<-\frac{c}{2}.
$$
In the last inequality we used \eqref{pani w szpileczkach}. This contradicts the positivity of $R$. \hfill $\square$

\subsection{Proof of Theorem \ref{tw d=2 independent coord.} } \label{proof_d2}

\noindent
{\bf Proof:} In view of Theorem \ref{TwNiez} the generating pairs $(G,Z)$ are such that 
\begin{gather}\label{równanie z Filipovica d=2}
	J_1(bG_1(x))+J_2(bG_2(x))=x \tilde{J}_{\mu}(b), \quad b,x\geq 0,
\end{gather}
where $\tilde{J}_{\mu}$ takes the form \eqref{pierwszy przyp Jmu} or \eqref{drugi przyp Jmu}. We deduce from \eqref{równanie z Filipovica d=2} the form of $G$ and characterize the noise $Z$. First let us consider the case when 
\begin{gather}\label{znikanie pochodnej}
\left(\frac{G_2(x)}{G_1(x)}\right)^\prime=0, \qquad x>0.
\end{gather}
Then 
$G(x)$ can be written in the form
\begin{gather*}
G(x)=g(x)\cdot\left(
\begin{array}{ccc}
 G_1\\
 G_2, 
\end{array}
\right), \quad x\geq 0,
\end{gather*}
with some function $g(x)\geq 0, x\geq 0$, and constants $G_1>0,G_2>0$. Equation \eqref{rownanie 2} 
amounts then to
\begin{align*}
dR(t)&=F(R(t))+g(R(t-)) \left(G_1dZ_1(t)+G_2 dZ_2(t)\right)\\[1ex]
&=F(R(t))+g(R(t-)) d\tilde{Z}(t), \quad t\geq 0,
\end{align*}
which is an equation driven by the one dimensional L\'evy process $\tilde{Z}(t):=G_1 Z_1(t)+G_2 Z_2(t)$. 
It follows that $\tilde{Z}$ is $\alpha_1$-stable with $\alpha_1\in(1,2]$ and that $g(x)=c_0x^{1/ \alpha_1}, c_0>0$. 
Notice that $Z^{G(x)}(t)=c_0 x^{\frac{1}{\alpha_1}}\tilde{Z}$, so $J_{Z^{G(x)}}(b)=c_{\alpha_1}(c_0 x^{\frac{1}{\alpha_1}}b)^{\alpha_1}=x c_0^{\alpha_1}c_{\alpha_1} b^{\alpha_1}$ and $c_0=(\frac{\eta_1}{c^{\alpha_1}})^{\frac{1}{\alpha_1}}$.
Hence \eqref{pierwszy przyp Jmu} holds and this proves $(Ia)$.

If \eqref{znikanie pochodnej} is not satisfied, then
\begin{gather}\label{niezerowanie pochodnej}
\left(\frac{G_2(x)}{G_1(x)}\right)^\prime\neq 0, \quad x\in (\underline{x},\bar{x}),
\end{gather}
for some interval $(\underline{x},\bar{x})\subset (0,+\infty)$. In the rest of the proof we consider this case and
prove $(Ib)$ and $(II)$.

$(Ib)$ From the equation
\begin{gather}\label{rrr}
J_1(bG_1(x))+J_2(bG_2(x))=x\eta_1 b^{\alpha_1}, \quad b\geq 0, \ x\geq 0,
\end{gather}
we explicitly determine unknown functions. Inserting $b/G_1(x)$ for $b$ yields
\begin{gather}\label{rowwww do eliminacji J1}
J_1(b)+J_2\left(b\frac{G_2(x)}{G_1(x)}\right)=\eta_1\frac{x}{G_1^{\alpha_1}(x)}b^{\alpha_1}, \quad b\geq 0, \quad x>0.
\end{gather}
Differentiation over $x$ yields
 $$
J_2^\prime\left(b\frac{G_2(x)}{G_1(x)}\right)\cdot b \left(\frac{G_2(x)}{G_1(x)}\right)^\prime
=\eta_1\left(\frac{x}{G_1^{\alpha_1}(x)}\right)^\prime b^{\alpha_1}, \quad b\geq 0 ,\quad  x>0.
 $$
Using \eqref{niezerowanie pochodnej} and dividing by $\left(\frac{G_2(x)}{G_1(x)}\right)^\prime$ leads to
 $$
J_2^\prime\left(b\frac{G_2(x)}{G_1(x)}\right)\cdot b=\eta_1 \frac{\left(\frac{x}{G^{\alpha_1}_1(x)}\right)^\prime}{\left(\frac{G_2(x)}{G_1(x)}\right)^\prime}\cdot b^{\alpha_1},\quad b\geq 0 ,\quad x\in(\underline{x},\bar{x}).
$$
By inserting $b\frac{G_1(x)}{G_2(x)}$ for $b$ one computes the derivative of $J_2$: 
$$
J_2^\prime(b)=\eta_1 \frac{\left(\frac{x}{G^{\alpha_1}_1(x)}\right)^\prime\left(\frac{G_1(x)}{G_2(x)}\right)^{\alpha_1-1}}{\left(\frac{G_2(x)}{G_1(x)}\right)^\prime}\cdot b^{\alpha_1-1},\quad b>0 , \quad x\in(\underline{x},\bar{x}).
$$
Fixing $x$ and integrating over $b$ provides
\begin{gather}\label{J_2 wylioczona}
J_2(b)=c_2 b^{\alpha_1}, \quad b>0,
\end{gather}
with some $c_2\geq 0$. Actually $c_2>0$ as $Z_2$ is of infinite variation and $J_2$ can not disappear.

By the symmetry of \eqref{rrr} the same conclusion holds for $J_1$, i.e.
\begin{gather}\label{J_1 wylioczona}
J_1(b)=c_1 b^{\alpha_1}, \quad b>0,
\end{gather}
with $c_1>0$. Using \eqref{J_2 wylioczona} and \eqref{J_1 wylioczona} in \eqref{rrr} gives us \eqref{Ib}. This proves $(Ib)$.

 $II)$ Solving the equation 
\begin{gather}\label{rrrrr}
J_1(bG_1(x))+J_2(bG_2(x))=x(\eta_1 b^{\alpha_1}+\eta_2 b^{\alpha_2}), \quad b,x\geq 0,
\end{gather}
in the same way as we solved \eqref{rrr} yields that 
 \begin{gather}\label{J1,J2 podwojne stabilne}
 J_1(b)=c_1 b^{\alpha_1}+c_2 b^{\alpha_2}, \quad  J_2(b)=d_1 b^{\alpha_1}+d_2 b^{\alpha_2}, \quad b\geq 0,
 \end{gather}
 with $c_1,c_2,d_1,d_2\geq 0$, $c_1+c_2>0, d_1+d_2>0$. From \eqref{rrrrr} and \eqref{J1,J2 podwojne stabilne} we can specify the following conditions for $G$:
  \begin{align}\label{aaaaaa}
c_1G_1^{\alpha_1}(x)+d_1G_2^{\alpha_1}(x)&=\eta_1 x,\\[1ex]\label{bbbbbb}
c_2G_1^{\alpha_2}(x)+d_2G_2^{\alpha_2}(x)&=\eta_2 x.
\end{align}
 We will show that $c_1>0, c_2=0, d_1=0, d_2>0$ by excluding the opposite cases.

 If $c_1>0,c_2>0$, one computes from \eqref{aaaaaa}-\eqref{bbbbbb} that
\begin{gather}\label{G2 na dwa sposoby}
G_1(x)=\left(\frac{1}{c_1}(\eta_1x-d_1G_2^{\alpha_1}(x))\right)^{\frac{1}{\alpha_1}}=\left(\frac{1}{c_2}(\eta_2x-d_2G_2^{\alpha_2}(x))\right)^{\frac{1}{\alpha_2}}, \quad x\geq 0.
\end{gather}
 This means that, for each $x\geq 0$, the value $G_2(x)$ is a solution of the following equation of the $y$-variable
\begin{gather}\label{rownanie z y}
\left(\frac{1}{c_1}(\eta_1x-d_1y^{\alpha_1})\right)^{\frac{1}{\alpha_1}}=\left(\frac{1}{c_2}(\eta_2x-d_2y^{\alpha_2})\right)^{\frac{1}{\alpha_2}}, 
\end{gather} 
with $y\in \left[0,\left(\frac{\gamma_1 x}{d_1}\right)^{\frac{1}{\alpha_1}}\wedge
\left(\frac{\gamma_2 x}{d_2}\right)^{\frac{1}{\alpha_2}}\right]$.
If $d_1=0$ or $d_2=0$ we compute $y=y(x)$ from \eqref{rownanie z y} and see that $d_1y^{\alpha_1}$ or $d_2y^{\alpha_2}$ must be negative either for $x$ sufficiently close to $0$ or $x$ sufficiently large. Now we need to exclude the case $d_1>0,d_2>0$. However,  in the case $c_1, c_2,d_1,d_2>0$ equation \eqref{rownanie z y} has no solutions because, for sufficiently large $x>0$, the left side of \eqref{rownanie z y} is strictly less than the right side. This inequality follows from Proposition \ref{prop o braku rozwiazan} proven below. 

So, we proved that $c_1\cdot c_2=0$ and similarly one proves that $d_1\cdot d_2=0$.
The case $c_1=0, c_2>0, d_1>0, d_2=0$ can be rejected because then $J_1$ would vary regularly with index $\alpha_2$
and $J_2$ with index $\alpha_1$, which is a contradiction. It follows that $c_1>0, c_2=0, d_1=0, d_2>0$ and in this case we obtain \eqref{postac g w drugim prrzypadku} from \eqref{aaaaaa} and \eqref{bbbbbb}. \hfill$\square$
\begin{prop}\label{prop o braku rozwiazan}
Let  $a,b,c,d>0$, $\gamma\in(0,1)$,  $2\geq \alpha_1>\alpha_2>1$. Then for sufficiently large $x>0$ the following inequalities are true
\begin{gather}\label{niernier pomocnicza}
\Big(ax-(bx-cz)^\gamma\Big)^{\frac{1}{\gamma}}-dz> 0, \qquad z\in \Big[0,\frac{b}{c}x\Big],
\end{gather}
\begin{gather}\label{niernier pomocnicza wlasciwa}
(bx-cy^{\alpha_1})^{\frac{1}{\alpha_1}}<(ax-dy^{\alpha_2})^{\frac{1}{\alpha_2}}, \quad y\in \Big[0,\Big(\frac{b}{c}x\Big)^{\frac{1}{\alpha_1}}\wedge \Big(\frac{a}{d}x\Big)^{\frac{1}{\alpha_2}}\Big].
\end{gather}
\end{prop}
{\bf Proof:} First we prove \eqref{niernier pomocnicza} and write it in the equivalent form
\begin{gather}\label{equivalentt}
ax\geq (dz)^{\gamma}+(bx-cz)^\gamma=:h(z).
\end{gather}
Since
$$
h^{\prime}(z)=\gamma \Big(d^\gamma z^{\gamma-1}-c(bx-cz)^{\gamma-1}\Big),
$$
$$
h^{\prime\prime}(z)=\gamma(\gamma-1)\Big(d^\gamma z^{\gamma-2}+c^2(bx-cz)^{\gamma-2}\Big)< 0, \quad z\in \Big[0,\frac{b}{c}x\Big],
$$
the function $h$ is concave and attains its maximum at point
$$
z_0:=\theta x:=\frac{b c^{\frac{1}{\gamma-1}}}{d^{\frac{\gamma}{\gamma-1}}+c^{\frac{\gamma}{\gamma-1}}}x \in \Big[0,\frac{b}{c}x\Big],
$$
which is a root of $h^\prime$. It follows that
\begin{align*}
h(z)\leq h(\theta x)&=(\theta x)^{\gamma}+(bx-c\theta x)^{\gamma}\\
&=(\theta^\gamma+(b-c\theta)^\gamma)x^{\gamma}< ax,
\end{align*}
provided that $x$ is sufficiently large and \eqref{niernier pomocnicza} follows. \eqref{niernier pomocnicza wlasciwa} follows from
\eqref{niernier pomocnicza} by setting $\gamma=\alpha_2/\alpha_1$, $z=y^{\alpha_1}$.
\hfill$\square$

\end{document}